\documentclass{aims}
\usepackage{amsmath}
\usepackage{lineno}
\usepackage{amsfonts}
\usepackage{algorithm}
\usepackage{algorithmicx}
\usepackage{empheq}
\usepackage{amsmath}
\usepackage{amsthm}
\usepackage{amsfonts}
\usepackage{subfig}
\usepackage{enumitem}
\usepackage{mathtools,cases}
\usepackage{mathtools,mathptmx}

\usepackage{paralist}
\usepackage{graphics} 
\usepackage{epsfig} 
\usepackage{graphicx}  \usepackage{epstopdf}
\usepackage[colorlinks=true]{hyperref}
\hypersetup{urlcolor=blue, citecolor=red}

\def\currentvolume{X}
 \def\currentissue{X}
  \def\currentyear{200X}
   \def\currentmonth{XX}
    \def\ppages{X--XX}
     \def\DOI{10.3934/xx.xx.xx.xx}

\newtheorem{theorem}{Theorem}[section]

\newtheorem{lemma}[theorem]{Lemma}

\newtheorem*{problem}{Problem}
\theoremstyle{definition}

\newtheorem{remark}{Remark}

\numberwithin{equation}{section}
\DeclareMathOperator{\sign}{sgn}

\title[1D CIP with
Single Measurement Data] 
      {Convexification for a 1D Hyperbolic Coefficient Inverse Problem with Single Measurement Data}

\author[Alexey V. Smirnov, Michael V. Klibanov and Loc H. Nguyen]{}


 \email{smirnov.phd@gmail.com}
 \email{mklibanv@uncc.edu}
 \email{loc.nguyen@uncc.edu}

\thanks{The work was supported by US Army Research Laboratory and US Army
Research Office grant W911NF-19-1-0044.}

\thanks{$^*$ Corresponding author: Michael Klibanov (mklibanv@uncc.edu)}

\begin{document}
\maketitle

\centerline{\scshape Alexey Smirnov}

\medskip

\centerline{\scshape Michael Klibanov$^*$ and Loc Nguyen}
\medskip
{\footnotesize
 \centerline{ Department of Mathematics and Statistics}
 \centerline{University of North Carolina Charlotte}
 \centerline{Charlotte, NC 28223, USA}
}

\bigskip

 \centerline{(Communicated by the associate editor name)}

\begin{abstract}
A version of the convexification numerical method for a Coefficient Inverse
Problem for a 1D hyperbolic PDE is presented. The data for this problem are
generated by a single measurement event. This method converges globally. The
most important element of the construction is the presence of the Carleman
Weight Function in a weighted Tikhonov-like functional. This functional is
strictly convex on a certain bounded set in a Hilbert space, and the
diameter of this set is an arbitrary positive number. The global convergence
of the gradient projection method is established. Computational results
demonstrate a good performance of the numerical method for noisy data.
\end{abstract}

\section{Introduction}

\label{sec:1}

We call a numerical method for a Coefficient Inverse Problem (CIP) \emph{%
globally convergent} if there exists a theorem claiming that this method
delivers at least one point in a sufficiently small neighborhood of the
exact solution without an assumption that the starting point of iterations
is located sufficiently close to that solution. We construct in this paper a
globally convergent numerical method for a CIP for a 1D hyperbolic PDE. This
CIP has a direct application in standoff imaging of dielectric constants of
explosive-like targets using experimentally collected data. Our numerical
method is a version of the so-called \emph{convexification} concept. Just as
in all previous publications about the convexification, which are cited
below, we work with the data resulting from a single measurement event.
Thus, our data depend on one variable.

The reason of our work on the convexification method is the well known fact that conventional Tikhonov least squares cost functionals for CIPs suffer from the phenomenon of multiple local minima and ravines, see, e.g. the work of Scales, Fischer and Smith \cite{scales1992global} for a convincing numerical example of this phenomenon. On the other hand, any version of the gradient method of the minimization of that functional stops at any local
minimum. Therefore, a numerical reconstruction technique, which is based on the minimization of that functional, is unreliable.

The convexification method for our particular CIP was not constructed in the past. Thus, we develop some new ideas here. The first new idea is to apply certain new changes of variables to the original problem to obtain a new Cauchy problem with the lateral Cauchy data for a quasilinear integro-differential equation with Volterra-like integrals in it. As soon as the solution of this problem is obtained, the target unknown coefficient can be computed by a simple backwards calculation. The second new idea is to obtain a new Carleman estimate for the principal part of the operator of that equation (Theorem 4.1). The Carleman Weight Function (CWF) in that estimate is also new. A surprising and newly observed property of that Carleman estimate is that a certain resulting integral, the one over an interval of a certain straight line, is non-negative. It is this property, which, in combination with the rest of that Carleman estimate, enables us to construct the \emph{key element} of the convexification, a globally strictly convex cost functional with the above mentioned CWF in it and then to prove the global convergence of our numerical method (Theorems 4.2-4.6). Since such a functional was not constructed for our CIP\ in the past, then both this construction and follow up Theorems 4.2-4.6 are also new.

Below $x\in \mathbb{R},\hspace{0.3em}t>0.$ Let the function $a(x)\in C^{1}(%
\mathbb{R})$ possesses the following properties: 
\begin{align}
& a(x)\geq 0\quad \mbox{for}\quad x\in (0,1),  \label{1.1} \\
& a(x)=0\quad \mbox{for}\quad x\notin (0,1).  \label{1.2}
\end{align}

\begin{problem}
(Forward Problem.) The forward problem we consider here is the problem of
the search of the fundamental solution $u(x,t)$ of the hyperbolic operator $%
\partial _{t}^{2}-\partial _{x}^{2}-a(x),$ with $a(x)$ satisfying (\ref{1.1}%
), (\ref{1.2}) i.e. 
\begin{equation}
\begin{dcases} \begin{aligned} &u_{tt}=u_{xx}+a(x) u, \quad (x,t) \in
\mathbb{R} \times (0,\infty ), \label{1.3}\\ &u(x,0) =0, \quad u_{t}(x,0)
=\delta (x), \end{aligned} \end{dcases}
\end{equation}%
where $\delta \left( x\right) $ is the Dirac function at $x=0.$
\end{problem}

\begin{problem}
(Coefficient Inverse Problem). \emph{\ Determine the coefficient }$a(x) $%
\emph{\ satisfying conditions (\ref{1.1}), (\ref{1.2}), assuming that the
following two functions }$f_{0}(t) ,f_{1}(t) $\emph{\ are given:}%
\begin{equation}
u(0,t) =f_{0}(t),\quad u_{x}(0,t) =f_{1}(t) ,\quad \forall t\in (0,T) ,
\label{1.4}
\end{equation}
\emph{where the number }$T>0$\emph{\ will be defined later.}
\end{problem}

It is the CIP (\ref{1.3}), (\ref{1.4}) for which we develop here the
convexification method. It is well known that, given (\ref{1.2}), functions $f_{0}(t),$ $f_{1}(t)$ for $t\in (0,2)$ (i.e. for $T=2$) uniquely determine the function $a(x)$ and also the Lipschitz stability estimate holds, see Theorem 2.6 Section 3 of Chapter 2 of \cite{romanov2018inverse} as well as Figure 1(B).

As to the Dirac function in the initial condition (\ref{1.3}),\emph{\ }this function is an
idealization of the reality of course. Therefore, its approximation is used in real world problems of physics. Nevertheless, the Dirac function is commonly used in many applied problems to model an ultra-short pulse, that penetrates deeply lossy materials and allows one to achieve very fine imaging resolution. An ultra-short pulse system is attractive for applications, due to its low power spectral density that results in negligible interference with other signals. There are various techniques to generate short pulses in the order of nanoseconds. In this regard, we refer to, e.g. an applied paper \cite{ahajjam2016compact}, where a short pulse is approximated via a narrow Gaussian.\ It is well known that such a function approximates the Dirac function in a certain sense. Another confirmation of the usefulness of the modeling via the Dirac function comes from \cite{klibanov2018kolesov}, where this function was successfully used to work with some experimental data via a version of the convexification method for a 1D CIP in the frequency domain.

To describe some applications of our CIP, we briefly consider here a
similar inverse problem for the 1D acoustic equation, 
\begin{equation}
\begin{dcases} \begin{aligned} &U_{tt}=c^{2}(y) U_{yy},\quad (y,t) \in
\mathbb{R} \times (0,\infty ) , \label{1.5}\\ &U(y,0) =0, \quad U_{t}(y,0)
=\delta (y). \end{aligned} \end{dcases}
\end{equation}%
where the sound speed $c(y)\in C^{3}(\mathbb{R})$ is such that $c(y)\geq
c_{0}=\mbox{const}>0$ and \newline
$c(y)=1$ for $y\in \left\{ \left( -\infty ,0\right) \cup (1,\infty )\right\} 
$. The coefficient inverse problem in this case consists of determining the
function $c(y)$ for $y\in (0,1),$ given functions $g_{0}(t)$ and $g_{1}(t),$%
\begin{equation}
U(0,t)=g_{0}(t),U_{y}(0,t)=g_{1}(t),\quad t\in (0,T^{\prime }),  \label{1.6}
\end{equation}%
where the number $T^{\prime }=T^{\prime }(T)$ depends on $T$ in (\ref{1.4}).

We start by applying a widely known change of variables, see e.g. \cite%
{romanov2018inverse}: 
\begin{equation*}
x\leftrightarrow y\quad \Rightarrow \quad x(y)=\int \displaylimits_{0}^{y}%
\frac{ds}{c(s)}
\end{equation*}%
Then $x(y)$ is the travel time of the acoustic signal from the point $%
\left\{ 0\right\} $ to the point $\left\{ y\right\} .$ Next, we introduce a
new function $V(x,t)=U(y(x),t)/S(x),$ where $S(x)=\sqrt{c(y(x))}.$ Then
problem (\ref{1.5})-(\ref{1.6}) becomes 
\begin{equation}
\begin{dcases} \begin{aligned} &V_{tt}=V_{xx}+p(x) V,\quad (x,t) \in
\mathbb{R} \times (0,\infty ) , \label{1.7}\\ &V(x,0) =0, \quad V_{t}(x,0)
=\delta (x),\\ &V(0,t) =g_{0}(t), \quad V_{x}(0,t)=g_{1}(t), \quad t\in
(0,T), \end{aligned} \end{dcases}
\end{equation}%
where 
\begin{equation*}
p(x)=\frac{S^{\prime \prime }(x)}{S(x)}-2\left[ \frac{S^{\prime }(x)}{S(x)}%
\right] ^{2}=\frac{1}{2}c^{\prime \prime }(y(x))c(y(x))-\frac{1}{4}\left[
c^{\prime }\left( y(x)\right) \right] ^{2}.
\end{equation*}

Equations (\ref{1.7}) look exactly as equations (\ref{1.3})-(\ref{1.4}).
Hence, we have reduced the CIP (\ref{1.5})-(\ref{1.6}) to our CIP (\ref{1.3}%
)-(\ref{1.4}). This justifies the applied aspect of our CIP. On the other
hand, due to the presence of the unknown coefficient $c(y)$ in the principal
part of the hyperbolic operator of (\ref{1.5}), the CIP (\ref{1.5})-(\ref%
{1.6}) is harder to work with than the CIP (\ref{1.3})-(\ref{1.4}).
Therefore, it makes sense, as the first step, to develop a numerical method
for the CIP (\ref{1.3})-(\ref{1.4}). Next, one might adapt that technique to
problem (\ref{1.5})-(\ref{1.6}). This first step is done in the current
paper.

The CIP (\ref{1.5})-(\ref{1.6}) has application in acoustics \cite{colton2019inverse}. Another quite interesting application is in inverse scattering of electromagnetic waves, in which case $c^{-2}(y)=\varepsilon
_{r}(y),$ where $\varepsilon _{r}(y)$ is the spatially distributed dielectric constant. Using the data, which were experimentally collected by the US Army Research Laboratory, it was demonstrated in \cite{Karch,klibanov2018kolesov,kuzhuget2012blind} that the 1D mathematical model, which is based on equation (\ref{1.5}), can be quite effectively used to image in the standoff mode dielectric constants of targets, which mimic explosives, such as, e.g. antipersonnel land mines and improvised explosive devices. In fact, the original data in \cite{Karch,klibanov2018kolesov,kuzhuget2012blind} were collected in the time domain. However, the mathematical apparatus of these references works only either with the Laplace transform \cite{Karch,kuzhuget2012blind} or with the
Fourier transform \cite{klibanov2018kolesov} with respect to $t$ of equation (\ref{1.5}). Unlike these, we hope that an appropriately modified technique of the current paper should help us in the future to work with those experimental
data directly in the time domain.

Of course, the knowledge of the dielectric constant alone is insufficient to differentiate between explosives and non-explosives. However, we believe that this knowledge might be used in the future as an ingredient, which would be an additional one to the currently existing features which are used in the classification procedures for such targets. So that this additional ingredient would decrease the current false alarm rate, see, e.g. page 33 of \cite{kuzhuget2012blind} for a similar conclusion. As to other globally convergent numerical methods for the 1D CIPs for the wave-like equations, we refer to works of Korpela, Lassas and Oksanen \cite{korpela2016regularization,korpela2018discrete}, where a CIP for equation (\ref{1.5}) is studied without the above change of variables. The data of \cite{korpela2016regularization,korpela2018discrete} depend on two variables since those are the Neumann-to-Dirichlet data. We also refer to the works of Kabanikhin with coauthors. First, this group has computationally implemented in the 1D case \cite{kabanikhin2013direct} the Gelfand-Krein-Levitan method (GKL) \cite{GL,Krein}. Next, they have extended the GKL method to the 2D case and studied that extension computationally, see, e.g. \cite{kabanikhin2013direct,kabanikhin2015fast,kabanikhin2015numerical}. In the original 1D version of GKL \cite{GL,Krein}, one reduces an analog of our CIP to a Fredholm-type linear integral equation of the second kind. The data for the CIP form the kernel of this equation. The solution of this equation provides one with the target unknown coefficient. In the 2D version of GKL, one obtains a system of coupled Fredholm-type linear integral equations of the second kind. The solution of this system allows one to calculate the unknown coefficient.

At the same time, it was demonstrated numerically in \cite{Karch} that while GKL works well for computationally simulated data in the 1D case, it fails to perform well for experimentally collected data. The latter is true at least for the experimental data of \cite{Karch}. These are the same experimental data as ones in \cite{klibanov2018kolesov,kuzhuget2012blind}.
This set of data is particularly important to us, since it is about the main application of our interest: imaging of dielectric constants of explosive-like targets. On the other hand, it was demonstrated in \cite{klibanov2018kolesov} that another 1D version of the convexification method performs well for the same experimental data. The version of \cite{klibanov2018kolesov} works with the data in the frequency domain, while the current paper works with the data in the time domain. We are not working with those experimental data in this paper, since such an effort would require a substantial investment of time from us, and we simply do not have this time at this moment. However, as stated above, in the future we indeed plan to apply the technique of the current paper to the experimental data of 
\cite{Karch,klibanov2018kolesov,kuzhuget2012blind}. Thus, we point out that while results of \cite{klibanov2018kolesov} show a good promise in this direction for the version of the convexification of the current paper, results of \cite{Karch} tell us that GKL is likely not applicable to those experimental data.

In the 2D case, the GKL uses overdetermined data \cite{kabanikhin2013direct,kabanikhin2015fast,kabanikhin2015numerical}. This means that the 2D version of GKL requires that the number $m=3$ of free variables in the data would exceed the number $n=2$ of free variables in the
unknown coefficient, i.e.$\ m>n$. On the other hand, in all publications about the convexification, which we cite below, so as in this one, the data are non overdetermined, i.e. $m=n.$ In particular, in this paper $m=n=1$.

Being motivated by the goal of avoiding the above mentioned phenomenon of multiple local minima and ravines of conventional least squares Tikhonov functionals, Klibanov with coauthors has been working on the convexification since 1995, see \cite{beilina2015globally,klibanov1995uniform,klibanov1997global,klibanov2016globally} for the initial works on this topic. The publication of Bakushinskii, Klibanov and Koshev \cite{bakushinskii2017carleman} has addressed some questions, which were important for the numerical implementation of the convexification. This has opened the door for some follow up publications about the convexification, including the current one, with a variety of computational results \cite{khoa2019convexification,Khoa2,klibanov2018kolesov,klibanov2019EIT,klibanov2019time,klibanov2018kolesov,klibanov2020parabolic}. We also refer to the works of Baudouin, De~Buhan and Ervedoza and Osses \cite{baudouin2017convergent,baudouin2020}, where a different version of the convexification is developed for two $n-$D CIPs ($n=1,2,...$) for the hyperbolic equations. Both versions of the convexification mentioned in this paragraph use the idea of the Bukhgeim-Klibanov method \cite{buchgeim1981uniqueness}.

As to the Bukhgeim-Klibanov method, it was originated in \cite{buchgeim1981uniqueness} with the only goal at that time (1981) of proofs of global uniqueness theorems for multidimensional CIPs with single measurement data. This method is based on Carleman estimates. The convexification extends the idea of \cite{buchgeim1981uniqueness} from the initial purely uniqueness topic to the more applied topic of numerical methods for CIPs. Many publications of many authors are devoted to the method of \cite{buchgeim1981uniqueness} being applied to a variety of CIPs, again with the goals of proofs of uniqueness and stability results for those CIPs. Since the current paper is not a survey of that technique, we now refer only to a few of such publications \cite{beilina2012approximate,klibanov1984inverse,klibanov1992inverse,klibanov2013carleman}%
.

All functions below are real valued ones. In Section 2 we derive a boundary
value problem for a quasilinear integro-differential equation. In Section 3
we describe the convexification method for solving this problem. We
formulate our theorems in Section 4. Their proofs are in Section 5.
Numerical results are presented in Section 6.

\section{Quasilinear Integro-Differential Equation}

\label{sec:2}

Let $H\left( x\right) $ be the Heaviside function centered at $x=0$. Problem
(\ref{1.3}) is equivalent to the following integral equation, see Section 3
of Chapter 2 of \cite{romanov2018inverse}:%
\begin{equation}
\begin{aligned} u\left( x,t\right) = \begin{dcases} \frac{1}{2}H\left(
t-\left\vert x\right\vert \right) +\frac{1}{2}\int\displaylimits_{D\left(
x,t\right) }a\left( \xi \right) u\left( \xi ,\tau \right) d\xi d\tau ,
&\text{ for }t>\left\vert x\right\vert,\\ 0,&\text{ for }0<t<\left\vert
x\right\vert.\end{dcases} \end{aligned}  \label{2.1}
\end{equation}%
\begin{equation}
D(x,t)=\left\{ (\xi ,\tau ):\left\vert \xi \right\vert <\tau <t-\left\vert
x-\xi \right\vert \right\} .  \label{2.2}
\end{equation}

\begin{figure}[tbp]
\begin{center}
\subfloat[$D(x,t) $]{\includegraphics[width
=.3\textwidth]{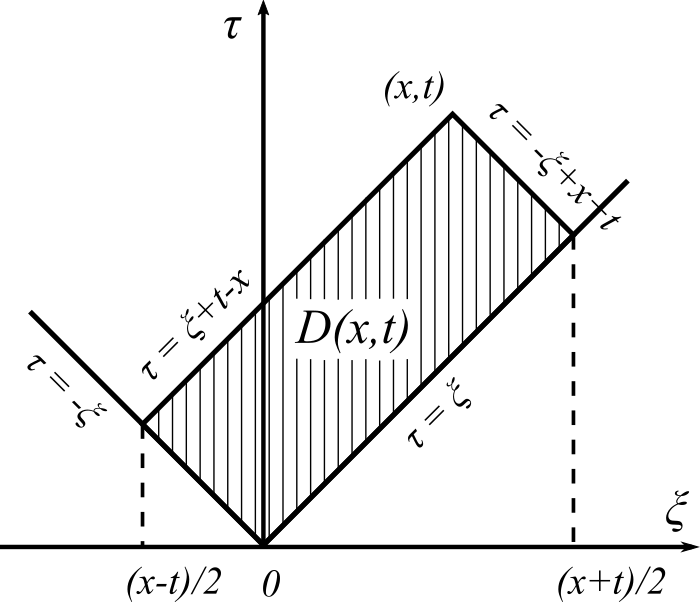}} \quad 
\subfloat[$D(0,t).$]{\includegraphics[width
=.3\textwidth]{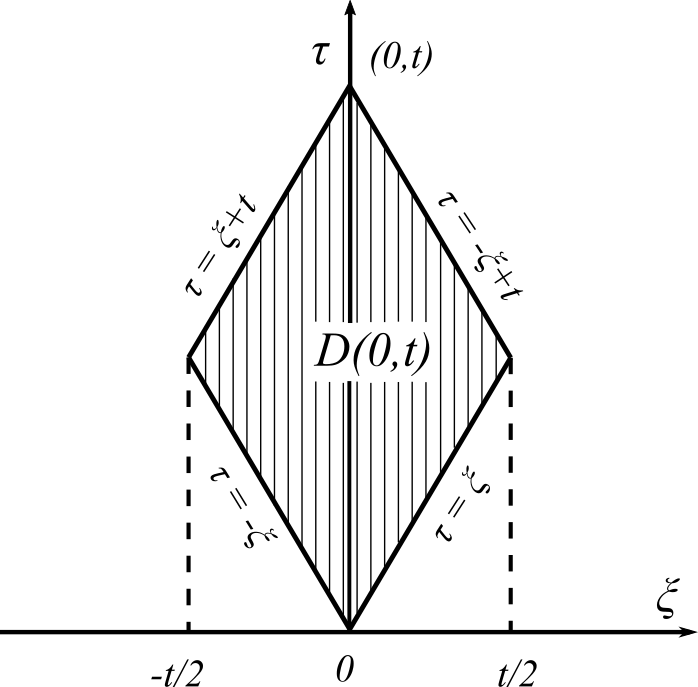}} \quad 
\subfloat[$Tr$ in (\ref{2.18})]{\includegraphics[width
=.255\textwidth]{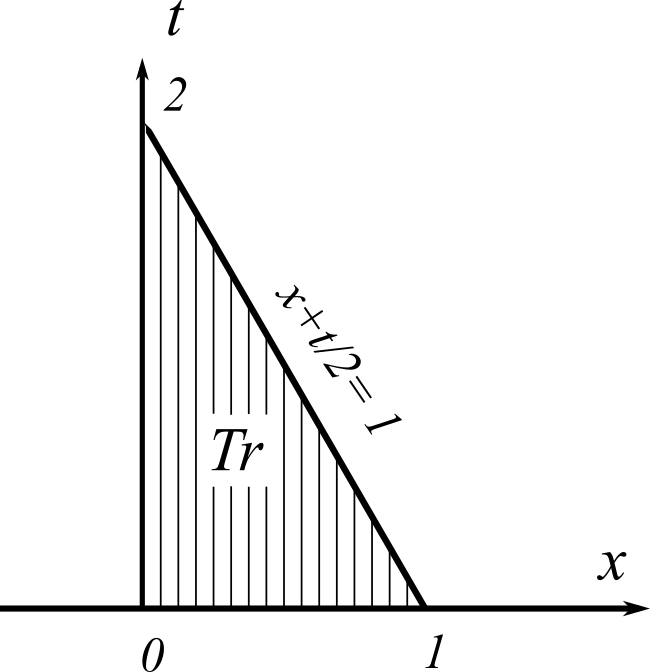}}
\end{center}
\caption{\emph{The rectangle} \textit{$D(x,t)=\left\{ (\protect\xi ,\protect%
\tau ):\left\vert \protect\xi \right\vert <\protect\tau <t-\left\vert x-%
\protect\xi \right\vert \right\} $ and the triangle $Tr$.}}
\label{fig1}
\end{figure}
It follows from (\ref{2.2}) and (\ref{1.2}) that the first line of (\ref{2.1}%
) can be rewritten as \cite{romanov2018inverse}: 
\begin{equation}
u(x,t)=\frac{1}{2}H(t-\left\vert x\right\vert )+\frac{1}{2}\int %
\displaylimits_{0}^{(x+t)/2}a(\xi )\int \displaylimits_{\left\vert \xi
\right\vert }^{t-\left\vert x-\xi \right\vert }u(\xi ,\tau )d\tau d\xi .
\label{2.3}
\end{equation}
see Figure 1. In fact, (\ref{2.3}) is a linear integral equation of the
Volterra type with respect to the function $u(x,t)$ \cite{romanov2018inverse}%
. This equation can be solved as:%
\begin{equation}
u_{0}(x,t)=\frac{1}{2}H(t-\left\vert x\right\vert ),\quad u_{n}(x,t)=\frac{1%
}{2}\int \displaylimits_{0}^{(x+t)/2}a(\xi )\int \displaylimits_{\left\vert
\xi \right\vert }^{t-\left\vert x-\xi \right\vert }u_{n-1}(\xi ,\tau )d\tau
d\xi   \label{2.4}
\end{equation}%
\begin{equation}
u(x,t)=\sum \displaylimits_{n=0}^{\infty }u_{n}(x,t),\quad \left\vert
u_{n}(x,t)\right\vert \leq \frac{(Mt)^{n}}{n!},\quad x\in (\alpha
_{1},\alpha _{2}),  \label{2.5}
\end{equation}%
for $n=1,2,\dots $ and for any finite interval $(\alpha _{1},\alpha
_{2})\subset \mathbb{R}$, where the number \newline
$M=M(\alpha _{1},\alpha _{2},\left\Vert a\right\Vert _{C[0,1]})>0$ depends
only on the listed parameters. Similar estimates can be obtained for
derivatives $\partial _{x}^{k}\partial _{t}^{s}u_{n}$ with $k+s\leq 3,$
except that in this case \newline
$M=M(\alpha _{1},\alpha _{2},\left\Vert a\right\Vert _{C^{1}[0,1]})>0.$ We
also note that since by (\ref{1.1}) $a(x)\geq 0,$ then (\ref{2.4})-(\ref{2.5}%
) imply that 
\begin{equation}
u(x,t)\geq \frac{1}{2}\text{ for }t\geq \left\vert x\right\vert .
\label{2.6}
\end{equation}

Thus, (\ref{2.1})-(\ref{2.6}) imply that the following lemma is valid \cite%
{romanov2018inverse}:

\begin{lemma}
\emph{There exists a unique solution } $u(x,t)$ \emph{of problem (\ref{2.1})
such that }$\left( u-u_{0}\right) (x,t)\in C\left\{ t\geq 0\right\} ,u\left(
x,t\right) \in C^{3}\left\{ (x,t)\hspace{0.3em}|\hspace{0.3em}t\geq
\left\vert x\right\vert \right\} $.\emph{\ Problem (\ref{2.1}) is equivalent
to the Cauchy problem (\ref{1.3})-(\ref{1.4}). Furthermore, }$%
\lim_{t\rightarrow \left\vert x\right\vert ^{+}}u(x,t)=1/2$ \emph{and} \emph{%
inequality (\ref{2.6}) holds.}
\end{lemma}

\subsection{Integro-differential equation}

\label{sec:2.1}

Consider the function $u(x,t)$ for $x>0$ above the characteristic cone $%
\left\{ t=\left\vert x\right\vert \right\} $ and change the variables as 
\begin{equation}
v(x,t)=u(x,t+x),\text{for }x,t>0.  \label{2.70}
\end{equation}%
Then (\ref{1.3}), (\ref{1.4}), (\ref{2.6}) and Lemma 2.1 imply that 
\begin{equation}
v_{xx}-2v_{xt}+a\left( x\right) v=0,\text{ for }x,t>0,  \label{2.8}
\end{equation}%
\begin{equation}
v\left( x,0\right) =\frac{1}{2},\text{ for }x>0,  \label{2.9}
\end{equation}%
\begin{equation}
v\left( 0,t\right) =f_{0}\left( t\right) ,v_{x}\left( 0,t\right)
=f_{0}^{\prime }\left( t\right) +f_{1}\left( t\right) .  \label{2.10}
\end{equation}%
In addition, (\ref{2.6}) and (\ref{2.70}) imply that 
\begin{equation}
v\left( x,t\right) \geq \frac{1}{2},\text{ for }x,t>0.  \label{2.11}
\end{equation}%
It follows from (\ref{2.11}) that we can consider the function 
\begin{equation}
q(x,t)=\ln v(x,t).  \label{2.110}
\end{equation}%
Using (\ref{2.8})-(\ref{2.10}), we obtain%
\begin{equation}
q_{xx}-2q_{xt}+q_{x}^{2}-2q_{x}q_{t}=-a\left( x\right) ,\text{ for }x,t>0,
\label{2.12}
\end{equation}%
\begin{equation}
q(x,t)=-\ln 2,  \label{2.13}
\end{equation}%
\begin{equation}
q\left( 0,t\right) =\ln f_{0}\left( t\right) ,\text{ }q_{x}\left( 0,t\right)
=\frac{f_{0}^{\prime }\left( t\right) +f_{1}\left( t\right) }{f_{0}\left(
t\right) }.  \label{2.14}
\end{equation}

Equation (\ref{2.12}) has two unknown functions, $q(x,t)$ and $a(x),$ which
is inconvenient. On the other hand, the function $a(x)$ is \textquotedblleft
isolated" \ in (\ref{2.12}) and it is independent on $t$. Therefore, we
follow the first step of the method of \cite{buchgeim1981uniqueness}. More
precisely, we differentiate both sides of equation (\ref{2.12}) with respect
to $t.$ Thus, we eliminate the unknown coefficient from this equation and
obtain an integro-differential equation this way.

Let 
\begin{equation}
w(x,t)=q_{t}(x,t).  \label{2.140}
\end{equation}%
Then (\ref{2.13}) and (\ref{2.140}) imply 
\begin{equation}
q(x,t)=\int \displaylimits_{0}^{t}w(x,\tau )d\tau -\ln 2.  \label{2.15}
\end{equation}%
Define the quasilinear integro-differential operator $L$ as 
\begin{equation}
L(w)=w_{xx}-2w_{xt}+2w_{x}\int \displaylimits_{0}^{t}w_{x}(x,\tau )d\tau
-2w_{x}w-2w_{t}\int \displaylimits_{0}^{t}w_{x}(x,\tau )d\tau .
\label{2.160}
\end{equation}%
Hence, (\ref{2.12})-(\ref{2.160}) imply 
\begin{equation}
L(w)=0,\text{ }\left( x,t\right) \in Tr,  \label{2.17}
\end{equation}%
\begin{equation}
w\left( 0,t\right) =p_{0}\left( t\right) ,\text{ }w_{x}\left( 0,t\right)
=p_{1}\left( t\right) ,  \label{2.170}
\end{equation}%
where%
\begin{equation}
p_{0}(t)=f_{0}^{\prime }(t)/f_{0}(t),\text{ }p_{1}(t)=\frac{d}{dt}%
[(f_{0}^{\prime }(t)+f_{1}(t))/f_{0}(t)]. \label{2.161}
\end{equation}
\vspace{0.5em}

As to the domain $Tr$ in (\ref{2.17}), it is clear that the change of
variables (\ref{2.70}) transforms the rectangle\emph{\ }$D(0,t)$ of Figure 1(B) in the triangle $Tr,$ see Figure 1(C),
\begin{equation}
Tr=\left\{ (x,t)\hspace{0.3em}:\hspace{0.3em}x,t>0,\hspace{0.3em}x+\frac{t}{2%
}<1\right\} .  \label{2.18}
\end{equation}%
Hence, we can uniquely determine the functions $w(x,t)$ and $q(x,t)$ only
for $(x,t)\in Tr.$

\subsection{Absorbing boundary conditions}

\label{sec:2.2}

\begin{lemma}
\emph{\ For every two numbers }$A\geq 1$ and $B> 0,$ \emph{the function }$%
u\left( x,t\right) $\emph{\ satisfies the absorbing boundary conditions: 
\begin{equation*}
\quad u_{x}(A,t)+u_{t}(A,t)=0,\text{ }u_{x}(-B,t)-u_{t}(-B,t)=0,\text{ }%
\forall t\in (0,T).
\end{equation*}%
}
\end{lemma}

\textbf{Proof}. Clearly the function $u_{0}\left( x,t\right) $ defined in (%
\ref{2.4}) satisfies these conditions. Denote $\widetilde{u}\left(
x,t\right) =u(x,t)-u_{0}(x,t).$ Differentiating (\ref{2.3}), we obtain%
\begin{equation}
\begin{aligned} &\widetilde{u}_{x}(x,t)=-\frac{1}{2}\int
\displaylimits_{0}^{\left( x+t\right) /2}\sign\left( x-\xi \right) a\left(
\xi \right) u\left( \xi ,t-\left\vert x-\xi \right\vert \right) d\xi , \\ &\widetilde{u}_{t}(x,t)=\frac{1}{2}\int
\displaylimits_{0}^{\left( x+t\right) /2}a\left( \xi \right) u\left( \xi
,t-\left\vert x-\xi \right\vert \right) d\xi . \end{aligned}  \label{100}
\end{equation}%
If $x\geq 1,$ then in (\ref{100}) $\sign\left( x-\xi \right) =1,$ since $%
a\left( \xi \right) =0$ for $\xi \geq 1.$ Next, if $x\leq 0,$ then in (\ref%
{100}) $\hspace{0.3em}\sign\left( x-\xi \right) =-1$ since $a\left( \xi
\right) =0$ for $\xi \leq 0.$ \ \ $\square $

\begin{remark}
\emph{Engquist and Majda have proposed to impose the absorbing boundary
conditions for the numerical simulations of the propagation of waves \cite%
{engquist1977absorbing}. Lemma 2.2 implies that, unlike \cite%
{engquist1977absorbing} , in the case of problem (\ref{1.3}), this condition
should not be imposed, since it holds automatically.}
\end{remark}

\begin{remark}
\emph{We impose the non-negativity condition (\ref{1.1}) on the unknown
coefficient }$a\left( x\right) $\emph{\ to ensure (\ref{2.6}). It is
inequality (\ref{2.6}), which allows us to consider the function }$%
q(x,t)=\ln v(x,t)$ in \emph{(\ref{2.110}): since (\ref{2.6}) guarantees (\ref%
{2.11}). Assumption (\ref{1.2}) is important for the validity of Lemma 2.2.
This lemma, in turn is quite helpful numerically for the solution of the
forward problem of data simulations as well as to ensure a good stability of
our inverse algorithm, see section 6. Finally, the smoothness condition }$%
a\in C^{1}\left( \mathbb{R}\right) $ \emph{ensures that the function }$q\in
C^{3}\left( x\geq 0,t\geq 0\right) :$\emph{\ see Lemma 2.1, (\ref{2.140})
and (\ref{2.160}). We point out that we are not looking for minimal
requirements imposed on }$a\left( x\right) .$
\end{remark}

Thus, (\ref{1.2}) and Lemma 2.2 imply that for any two numbers $A\geq 1,B>0$%
\begin{equation}
u_{tt}=u_{xx}+a\left( x\right) u,\text{ }\left( x,t\right) \in \left(
-B,A\right) \times \left( 0,\infty \right) ,  \label{2.180}
\end{equation}%
\begin{equation}
u\left( x,0\right) =0,u_{t}\left( x,0\right) =\delta \left( x\right) ,
\label{2.181}
\end{equation}%
\begin{equation}
u_{x}\left( -B,t\right) -u_{t}\left( -B,t\right) =0,\text{ }u_{x}\left(
A,t\right) +u_{t}\left( A,t\right) =0.  \label{2.182}
\end{equation}

\subsection{Reconstruction of the unknown coefficient}

\label{sec:2.3}

It follows from (\ref{2.12}), (\ref{2.13}) and (\ref{2.140}) that\emph{\ }%
\begin{equation}
a(x)=2w_{x}(x,0).  \label{2.19}
\end{equation}%
Hence, we focus below on the numerical solution of the boundary value
problem (\ref{2.17}), (\ref{2.161}).

\section{Convexification}

\label{sec:3}

\subsection{Convexification in brief}

\label{sec:3.1}

Given a CIP, the first step of the convexification follows the first step of \cite{buchgeim1981uniqueness}, in which the unknown coefficient is eliminated from the PDE via the differentiation with respect to such a parameter from which that coefficient does not depend. In particular, in our
case, we have replaced equation (\ref{2.12}), which contains the unknown coefficient $a(x),$ with a quasilinear integro-differential equation (\ref{2.17}), which does not contain that coefficient. Next, one should solve the
corresponding boundary value problem, which is similar with the problem (\ref{2.17}), (\ref{2.170}). To solve that boundary value problem, a weighted Tikhonov-like functional $J_{\lambda }$ is constructed, where $\lambda \geq 1$ is a parameter. The weight is the Carleman Weight Function (CWF), which is involved in the Carleman estimate for the principal part of the operator of that integro-differential equation. In our case, that principal part is the operator $\partial _{x}^{2}-2\partial _{x}\partial _{t},$ see (\ref%
{2.160}) and (\ref{2.17}).

The above mentioned functional is minimized on a convex bounded set with the
diameter $2d,$ where $d>0$ is an arbitrary number. This set is a part of a
Hilbert space $H^{k}.$ In our case, $k=3$. The key theorem is that one can
choose a sufficiently large value $\widetilde{\lambda }(d)\geq 1$ of the
parameter $\lambda $ such that the functional $J_{\lambda }$ is strictly
convex on that set for all $\lambda \geq \widetilde{\lambda }.$ Next, one
proves that, for these values of $\lambda ,$ the gradient projection method
of the minimization of the functional $J_{\lambda }$ converges to the
correct solution of that CIP starting from an arbitrary point of the above
mentioned set, as long as the level of the noise in the data tends to zero.
Given that the diameter $2d$ of that set is an arbitrary number and that the
starting point is also an arbitrary one, this is the \emph{global convergence%
}, by the definition of the first sentence of Introduction.

It is worth to note that even though the theory says that the parameter $%
\lambda $ should be sufficiently large, our rich computational experience
tells us that computations are far less pessimistic than the theory is. More
precisely, in all our numerically oriented publications on the
convexification, including the current one, accurate numerical results are
obtained for $\lambda \in \lbrack 1,3]$, see \cite{bakushinskii2017carleman,khoa2019convexification,klibanov2019convexification,klibanov2019time,klibanov2019EIT,klibanov2008new,klibanov2018kolesov}.

\subsection{The Tikhonov-like functional with the Carleman Weight Function in it}

\label{sec:3.2}

We construct this functional to solve problem (\ref{2.17}), (\ref{2.170}). Everywhere below $\alpha \in (0,1/2).$ Our CWF has the form: 
\begin{equation}
\varphi _{\lambda }(x,t)=\exp \left( -2\lambda (x+\alpha t)\right) ,
\label{3.1}
\end{equation}%
where $\lambda \geq 1$ is a parameter, see Theorem 4.1 in section 4 for the
Carleman estimate with this CWF. Even though we can find the function $%
w(x,t) $ only in the triangle $Tr$ in (\ref{2.18}), it is convenient for our
numerical study to work with the rectangle $R,$ 
\begin{equation}
R=(0,1)\times (0,T),\quad T\geq 2.  \label{3.2}
\end{equation}%
Using (\ref{2.70}), (\ref{2.110}), (\ref{2.140}) and the absorbing boundary
condition (\ref{2.182}) for $A=1,$ we obtain 
\begin{equation}
w_{x}\left( 1,t\right) =0.  \label{6.4}
\end{equation}%
Let $d>0$ be an arbitrary number. Define the set $B(d,p_{0},p_{1})$ as 
\begin{equation}
\begin{aligned} &B(d,p_{0},p_{1})=\\ &\left\{ w\in
H^{3}(R):w(0,t)=p_{0}(t),\hspace{0.3em}w_{x}(0,t)=p_{1}(t),%
\hspace{0.3em}w_{x}\left( 1,t\right) =0,\hspace{0.3em}\left\Vert
w\right\Vert _{H^{3}(R)}<d\right\} . \label{3.3}
\end{aligned}
\end{equation}%
Let $\beta \in (0,1)$ be the regularization parameter and $L(w)$ be the
operator defined in (\ref{2.160}). Our weighted Tikhonov-like functional is:%
\begin{equation}
J_{\lambda ,\beta }(w)=\int \displaylimits_{R}[L(w)]^{2}\varphi _{\lambda
}dxdt+\beta \left\Vert w\right\Vert _{H^{3}(R)}^{2}.  \label{3.4}
\end{equation}%
\textbf{Minimization Problem}. \emph{Minimize the functional }$J_{\lambda
,\beta }(w)$\emph{\ on the set }$B(d,p_{0},p_{1}).$

\subsection{Estimating an integral}

\label{sec:3.3}

We use Lemma 3.1 in the proof of Theorem 4.2 (section 4). The presence of the multiplier $1/\lambda ^{2}$ in the right hand side of (\ref{3.5}) is new since the CWF is new here. Indeed, while in (\ref{3.1}) $t$ is used, usually one uses $t^{2}$ in CWFs for similar problems, see e.g. \cite{beilina2012approximate,klibanov2013carleman}. The latter implies that the term $1/\lambda $ rather than $1/\lambda ^{2}$ is present in an analogous estimate of Lemma 1.10.3 of \cite{beilina2012approximate} and of Lemma 3.1 of \cite{klibanov2013carleman}. Since these and similar lemmata are usually used in the Bukhgeim-Klibanov method and since any Carleman estimate requires that its parameter $\lambda \geq 1$ be sufficiently large, then the estimate of Lemma 3.1 is stronger than the one of \cite{beilina2012approximate,klibanov2013carleman}. The proof of this estimate is also different from the one of \cite{beilina2012approximate,klibanov2013carleman}. Even though we use an arbitrary $\alpha >0$ in Lemma 3.1, still everywhere after this lemma $\alpha \in \left( 0,1/2\right) :$ just as above$.$

\begin{lemma}
For any two numbers $\lambda ,\alpha >0$ and \emph{for any function }$g\in
L^{2}(R)$\emph{\ the following estimate is valid:}%
\begin{equation}
\int \displaylimits_{R}\left( \int \displaylimits_{0}^{t}g(x,\tau )d\tau
\right) ^{2}\varphi _{\lambda }dxdt\leq \frac{1}{\lambda ^{2}\alpha ^{2}}%
\int \displaylimits_{R}g^{2}\varphi _{\lambda }dxdt.  \label{3.5}
\end{equation}
\end{lemma}

\textbf{Proof}. Using (\ref{3.1}), integration by parts and the
Cauchy-Schwarz inequality, we obtain%
\begin{align*}
&I =\int \displaylimits_{R}\left( \int \displaylimits_{0}^{t}g(x,\tau )d\tau
\right) ^{2}\varphi _{\lambda }dxdt=\int \displaylimits_{0}^{1}e^{-2\lambda
x}\int \displaylimits_{0}^{T}e^{-2\lambda \alpha t}\left( \int \displaylimits%
_{0}^{t}g(x,\tau )d\tau \right) ^{2}dtdx = \\
& \int \displaylimits_{0}^{1}e^{-2\lambda x}\int \displaylimits_{0}^{T}\frac{%
d}{dt}\left( -\frac{e^{-2\lambda \alpha t}}{2\lambda \alpha }\right) \left(
\int \displaylimits_{0}^{t}g(x,\tau )d\tau \right) ^{2}dtdx = \\
&-\int \displaylimits_{0}^{1}e^{-2\lambda x}\frac{e^{-2\lambda \alpha T}}{%
2\lambda \alpha }\left( \int \displaylimits_{0}^{T}g(x,\tau )d\tau \right)
^{2}dx \\
&+\frac{1}{\lambda \alpha }\int \displaylimits_{R}e^{-2\lambda
x}e^{-2\lambda \alpha t}g(x,t)\left( \int \displaylimits_{0}^{t}g(x,\tau
)d\tau \right) dtdx \leq \\
& \frac{1}{\lambda \alpha }\left[ \int \displaylimits_{R}g^{2}\varphi
_{\lambda }dxdt\right] ^{1/2}\left[ \int \displaylimits_{R}\left( \int %
\displaylimits_{0}^{t}g(x,\tau )d\tau \right) ^{2}\varphi _{\lambda }dxdt%
\right] ^{1/2}.
\end{align*}
Here, we have used the fact that the term in the third line of the above is negative. Hence, we have obtained that 
\begin{equation}
I\leq \frac{1}{\lambda \alpha }\left( \int \displaylimits_{R}g^{2}\varphi
_{\lambda }dxdt\right) ^{1/2}\sqrt{I}.  \label{3.6}
\end{equation}

Dividing both sides of (\ref{3.6}) by $\sqrt{I}$ and squaring both sides of the resulting inequality, we obtain (\ref{3.5}). \ $\square $

\section{Theorems}

\label{sec:4}

Introduce the subspaces $H_{0}^{2}(R)\subset H^{2}(R)$ and $%
H_{0}^{3}(R)\subset H^{3}(R),$ 
\begin{equation*}
H_{0}^{2}(R)=\left\{ u\in H^{2}(R):u(0,t)=u_{x}(0,t)\right\} ,\quad
H_{0}^{3}(R)=H^{3}(R)\cap H_{0}^{2}(R).
\end{equation*}

\begin{theorem}
(Carleman estimate). \emph{There exist constants }$C=C(\alpha )>0$\emph{\
and }$\lambda _{0}=\lambda _{0}(\alpha )\geq 1$\emph{\ depending only on }$%
\alpha $\emph{\ such that for all functions }$u\in H_{0}^{2}(R)$\emph{\ and
for all }$\lambda \geq \lambda _{0}$\emph{\ the following Carleman estimate
is valid:} 
\begin{equation}
\begin{aligned} &\int\displaylimits_{R}( u_{xx}-2u_{xt}) ^{2}\varphi
_{\lambda }dxdt\geq C\lambda \int\displaylimits_{R}( u_{x}^{2}+u_{t}^{2})
\varphi _{\lambda }dxdt+C\lambda ^{3}\int\displaylimits_{R}u^{2}\varphi
_{\lambda }dxdt \\ &+C\lambda \int\displaylimits_{0}^{1}u_{x}^{2}( x,0)
e^{-2\lambda x}dx+C\lambda ^{3}\int\displaylimits_{0}^{1}u^{2}( x,0)
e^{-2\lambda x}dx -C\lambda e^{-2\lambda \alpha
T}\int\displaylimits_{0}^{1}u_{x}^{2}( x,T) dx\\ &-C\lambda ^{3}e^{-2\lambda
\alpha T}\int\displaylimits_{0}^{1}u^{2}( x,T) dx.\label{4.1} \end{aligned}
\end{equation}
\end{theorem}

\begin{remark}
\emph{This Carleman estimate is new. The positivity of the first two terms
in the second line of (\ref{4.1}) is surprising. Indeed, in Carleman
estimates, usually one cannot ensure signs of integrals over hypersurfaces.
In particular, using (\ref{2.19}), it is shown below that the positivity of
these two terms is quite helpful in the reconstruction of the unknown
coefficient }$a(x).$ 
\end{remark}

Choose an arbitrary number $\varepsilon \in ( 0,2\alpha ) .$ Consider the
triangle $Tr_{\alpha ,\varepsilon }$%
\begin{equation}
Tr_{\alpha ,\varepsilon }=\left\{ ( x,t) :x+\alpha t<2\alpha -\varepsilon
;\quad x,t>0\right\} \subset Tr  \label{4.2}
\end{equation}

\begin{theorem}
(global strict convexity). \emph{For an arbitrary number }$d>0,$\emph{\ let }%
$B(d,p_{0},p_{1})\subset H^{3}(R)$\emph{\ be the set defined in (\ref{3.3}).
For any }$\lambda ,\beta >0$\emph{\ and for any }$w\in \overline{%
B(d,p_{0},p_{1})}$\emph{\ the functional }$J_{\lambda ,\beta }(w)$\emph{\ in
(\ref{3.4}) has the Fr\'{e}chet derivative }$J_{\lambda ,\beta }^{\prime
}(w)\in H_{0}^{3}(R).$\emph{\ Let }$\lambda _{0}=\lambda _{0}(\alpha )\geq 1$%
\emph{\ be the number of Theorem 4.1. Then there exist a sufficiently large
number }$\lambda _{1}=\lambda _{1}(\alpha ,\varepsilon ,d)\geq $\emph{\ }$%
\lambda _{0}$\emph{\ and a number }$C_{1}=C_{1}(\alpha ,\varepsilon ,d)>0$%
\emph{, both depending only on listed parameters, such that for all }$%
\lambda \geq \lambda _{1}$\emph{\ and for all }$\beta \in \lbrack
2e^{-\lambda \alpha T},1),$\emph{\ functional (\ref{3.4}) is strictly convex
on the set }$\overline{B(d,p_{0},p_{1})}$\emph{.\ More precisely, the
following inequality holds:}%
\begin{equation}
\begin{aligned} J_{\lambda ,\beta }\left( w_{2}\right) -J_{\lambda ,\beta
}\left( w_{1}\right) -J_{\lambda ,\beta }^{\prime }\left( w_{1}\right)
\left( w_{2}-w_{1}\right) \geq C_{1}e^{-2\lambda \left( 2\alpha -\varepsilon
\right) }\left\Vert w_{2}-w_{1}\right\Vert _{H^{1}\left( Tr_{\alpha
,\varepsilon }\right) }^{2} \label{4.3} \\ +C_{1}e^{-2\lambda \left( 2\alpha
-\varepsilon \right) }\left\Vert w_{2}\left( x,0\right) -w_{1}\left(
x,0\right) \right\Vert _{H^{1}\left( 0,2\alpha -\varepsilon \right)
}^{2}+\frac{\beta }{2}\left\Vert w_{2}-w_{1}\right\Vert _{H^{3}\left(
R\right) }^{2}, \\ \forall w_{1},w_{2}\in \overline{B\left(
d,p_{0},p_{1}\right) }, \hspace{0.3em} \forall \lambda \geq \lambda _{1}.
\end{aligned}
\end{equation}
\end{theorem}

\begin{remark}
\emph{Below }$C_{1}=C_{1}(\alpha ,\varepsilon ,d)>0$\emph{\ denotes
different numbers depending only on listed parameters.} \emph{It follows
from Lemma 3 on page 9 of the book of Polyak \cite{polyak1987introduction}
that (\ref{4.3}) guarantees the strict convexity of the functional }$%
J_{\lambda ,\beta }$ \emph{on the set} $\overline{B\left(
d,p_{0},p_{1}\right) }.$
\end{remark}

\begin{theorem}
\emph{Let parameters }$\lambda _{1},\lambda ,\beta $\emph{\ be the same as
in Theorem 4.2. Then there exists a unique minimizer }$w_{\min ,\lambda
,\beta }\in \overline{B(d,p_{0},p_{1})}$\emph{\ of the functional }$%
J_{\lambda ,\beta }(w)$\emph{\ on the set }$\overline{B(d,p_{0},p_{1})}.$%
\emph{\ Furthermore, the following inequality holds}%
\begin{equation}
J_{\lambda ,\beta }^{\prime }(w_{\min ,\lambda ,\beta })(w-w_{\min ,\lambda
,\beta })\geq 0,\quad \forall w\in \overline{B(d,p_{0},p_{1})}.  \label{4.4}
\end{equation}
\end{theorem}

To estimate the reconstruction accuracy as well as to introduce the gradient
projection method, we need to obtain zero Dirichlet and Neumann boundary
conditions at $\left\{ x=0\right\} .$ Also, we need to introduce noise in
the data and to consider an exact, noiseless solution. By one of the
concepts of the regularization theory, we assume that there exists an exact
solution $a^{\ast }(x)\in C^{1}(\mathbb{R})$ of the CIP (\ref{1.3})-(\ref%
{1.4}) with the noiseless data \cite%
{beilina2012approximate,tikhonov2013numerical}, and this function satisfies
conditions (\ref{1.1}), (\ref{1.2}). Let $w^{\ast }$ be the function $w$
which corresponds to $a^{\ast }(x)$. We assume that $w^{\ast }\in
B(d,p_{0}^{\ast },p_{1}^{\ast }),$ where $p_{0}^{\ast },p_{1}^{\ast }$ are
the noiseless data $p_{0},p_{1}.$ Let $\xi \in (0,1)$ be the level of noise
in the data. Obviously there exists a function $G^{\ast }\in B(d,p_{0}^{\ast
},p_{1}^{\ast }).$ Suppose that there exists a function $G\in
B(d,p_{0},p_{1})$ such that 
\begin{equation}
\left\Vert G-G^{\ast }\right\Vert _{H^{3}(R)}<\xi .  \label{4.5}
\end{equation}%
Denote $W^{\ast }=w^{\ast }-G^{\ast }$ and $W=w-G,$ $\forall w\in
B(d,p_{0},p_{1}),$ 
\begin{equation*}
B_{0}(D)=\left\{ U\in H_{0}^{3}(R):\left\Vert U\right\Vert
_{H^{3}(R)}<D\right\} ,\quad \forall D>0.
\end{equation*}%
Then (\ref{3.3}) and the triangle inequality imply that%
\begin{align}
& W^{\ast }\in B_{0}(2d),\quad W\in B_{0}(2d),\quad \forall w\in
B(d,p_{0},p_{1}),  \label{4.8} \\
& W+G\in B(3d,p_{0},p_{1}),\text{ }\forall W\in B_{0}(2d).  \label{4.9}
\end{align}%
Denote 
\begin{equation*}
I_{\lambda ,\beta }(W)=J_{\lambda ,\beta }(W+G),\hspace{0.3em}\forall W\in
B_{0}(2d).
\end{equation*}%
\ 

\begin{theorem}
\emph{The Fr\'{e}chet derivative }$I_{\lambda ,\beta }^{\prime }(W)\in
H_{0}^{3}(R)$\emph{\ of the functional }$I_{\lambda ,\beta }(W)$\emph{\
exists for every point }$W\in \overline{B_{0}(2d)}$\emph{\ and for all }$%
\lambda ,\beta >0.$\emph{\ Let }$\lambda _{1}=\lambda _{1}(\alpha
,\varepsilon ,d)$\emph{\ be the number of Theorem 4.2.\ Denote }$\lambda
_{2}=\lambda _{1}(\alpha ,\varepsilon ,3d)\geq \lambda _{1}.$\emph{\ Let }$%
\lambda \geq \lambda _{2}$\emph{\ and also let }$\beta \in \lbrack
2e^{-\lambda \alpha T},1).$\emph{\ Then the functional }$I_{\lambda ,\beta
}(W)$\emph{\ is strictly convex on the ball }$\overline{B_{0}(2d)}\subset
H_{0}^{3}(R).$\emph{\ More precisely, the following estimate holds:}%
\begin{equation}
\begin{aligned} I_{\lambda ,\beta }\left( W_{2}\right) -I_{\lambda ,\beta
}\left( W_{1}\right) -I_{\lambda ,\beta }^{\prime }\left( W_{1}\right)
\left( W_{2}-W_{1}\right) \geq C_{1}e^{-2\lambda \left( 2\alpha -\varepsilon
\right) }\left\Vert W_{2}-W_{1}\right\Vert _{H^{1}\left( Tr_{\alpha
,\varepsilon }\right) }^{2} \label{4.10} \\ +C_{1}e^{-2\lambda \left(
2\alpha -\varepsilon \right) }\left\Vert W_{2}\left( x,0\right) -W_{1}\left(
x,0\right) \right\Vert _{H^{1}\left( 0,2\alpha -\varepsilon \right)
}^{2}+\frac{\beta }{2}\left\Vert W_{2}-W_{1}\right\Vert _{H^{3}\left(
R\right) }^{2}, \\ \forall w_{1},w_{2}\in \overline{B_{0}\left( 2d\right) },
\hspace{0.3em}\forall \lambda \geq \lambda _{2}. \end{aligned}
\end{equation}

\emph{Furthermore, there exists a unique minimized }$W_{\min ,\lambda ,\beta
}$\emph{\ }$\in \overline{B_{0}(2d)}$\emph{\ of the functional }$I_{\lambda
,\beta }(W)$\emph{\ and the following inequality holds}%
\begin{equation}
I_{\lambda ,\beta }^{\prime }(W_{\min ,\lambda ,\beta })(W-W_{\min ,\lambda
,\beta })\geq 0,\quad \forall W\in \overline{B_{0}(2d)}.  \label{4.11}
\end{equation}
\end{theorem}

\vspace{0.3em}

\begin{theorem}
(the accuracy of the minimizer).\emph{\ Let the number }$T\geq 4.$\emph{\
Denote }%
\begin{equation}
\sigma =\frac{\alpha (T-4)+\varepsilon }{2(2\alpha -\varepsilon )},\quad
\rho =\frac{1}{2}\min (\sigma ,1)\in \left( 0,\frac{1}{2}\right) .
\label{4.12}
\end{equation}%
\emph{Choose a number }$\xi _{0}\in (0,1)$\emph{\ so small that }$\ln \xi
_{0}^{-1/(2(2\alpha -\varepsilon ))}\geq \lambda _{2},$\emph{\ where }$%
\lambda _{2}$\emph{\ is the number of} \emph{Theorem 4.4.} \emph{Let the
level of noise in the data }$\xi \in (0,\xi _{0}).$\emph{\ Choose the
parameters }$\lambda =\lambda (\xi )$ and $\beta =\beta (\xi )$ \emph{as} 
\emph{\ }%
\begin{equation}
\lambda =\lambda (\xi )=\ln \xi ^{-1/(2(2\alpha -\varepsilon ))}>\lambda
_{2},\quad \beta =\beta (\xi )=2e^{-\lambda \alpha T}=2\xi ^{(\alpha
T)/(2(2\alpha -\varepsilon ))}  \label{4.121}
\end{equation}%
\emph{\ (see Theorem 4.2 for }$\beta $\emph{). Then the following accuracy
estimates are valid:}%
\begin{equation}
\left\Vert w_{\min ,\lambda ,\beta }-w^{\ast }\right\Vert _{H^{1}(Tr_{\alpha
,\varepsilon })}\leq C_{1}\xi ^{\rho },\quad \left\Vert a_{\min ,\lambda
,\beta }-a^{\ast }\right\Vert _{L^{2}(0,2\alpha -\varepsilon )}\leq C_{1}\xi
^{\rho },  \label{4.14}
\end{equation}%
\emph{where }$w_{\min ,\lambda ,\beta }=(W_{\min ,\lambda ,\beta }+G)\in 
\overline{B(3d,p_{0},p_{1})}.$ \emph{Here,} $W_{\min ,\lambda ,\beta }\in 
\overline{B_{0}(2d)}$\emph{\ is the minimizer, which is found in Theorem
4.4, and }$a_{\min ,\lambda ,\beta }(x)=2\partial _{x}[w_{\min ,\lambda
,\beta }(x,0)],$ \emph{as in (\ref{2.19}). }
\end{theorem}

We now construct the gradient projection method of the minimization of the
functional $I_{\lambda ,\beta }(W)$ on the closed ball $\overline{B_{0}(2d)}%
\subset H_{0}^{3}(R).$ Let $P_{B_{0}}:H_{0}^{3}(R)\rightarrow \overline{%
B_{0}(2d)}$ be the orthogonal projection operator. Let $W_{0}\in B_{0}(2d)$
be an arbitrary point and the number $\gamma \in (0,1).$ The sequence of the
gradient projection method is \cite{bakushinskii2017carleman}:%
\begin{equation}
W_{n}=P_{B_{0}}(W_{n-1}-\gamma I_{\lambda ,\beta }^{\prime }(W_{n-1})),\quad
n=1,2,...  \label{4.15}
\end{equation}

\begin{theorem}
(the global convergence of the gradient projection method). \newline
\emph{Let }$\lambda _{2}=\lambda _{1}(\alpha ,\varepsilon ,3d)\geq \lambda
_{1},$\emph{\ where }$\lambda _{1}\geq 1$\emph{\ is the number of Theorem
4.2. Let }$\ $\emph{the numbers }$T$\emph{,}$\rho ,\xi _{0},\xi \in (0,\xi
_{0}),\lambda (\xi )$\emph{\ and }$\beta (\xi )$\emph{\ be the same as in
Theorem 4.5. Let }$W_{\min ,\lambda ,\beta }\in \overline{B_{0}(2d)}$\emph{\
be the unique minimizer of the functional }$I_{\lambda ,\beta }(W),$\emph{\
as in Theorem 4.4. Also, as in Theorem 4.4, denote }$w_{\min ,\lambda ,\beta
}=(W_{\min ,\lambda ,\beta }+G)\in \overline{B(3d,p_{0},p_{1})}$\emph{\ and
let }$w_{n}=(W_{n}+G)\in \overline{B(3d,p_{0},p_{1})},$\emph{\ where }$%
n=0,1,....$ \emph{Also, let }$a_{\min ,\lambda ,\beta }(x)$\emph{\ and }$%
a_{n}(x)$\emph{\ be the approximations of the coefficient }$a^{\ast }(x),$%
\emph{\ which are found from the functions }$w_{\min ,\lambda ,\beta }$\emph{%
\ and }$w_{n}$\emph{\ respectively via (\ref{2.19}). Then there exists a
number }$\gamma _{0}=\gamma _{0}(\alpha ,\varepsilon ,d,\xi )\in (0,1)$\emph{%
\ depending only on listed parameters such that for any }$\gamma \in
(0,\gamma _{0})$\emph{\ there exists a number }$\theta =\theta (\gamma )\in
(0,1)$\emph{\ such that the following convergence rates hold:}%
\begin{align}
& \left\Vert w_{\min ,\lambda ,\beta }-w_{n}\right\Vert _{H^{3}(R)}\leq
\theta ^{n}\left\Vert w_{\min ,\lambda ,\beta }-w_{0}\right\Vert
_{H^{3}(R)},\quad n=1,2,...,  \label{4.16} \\
& \left\Vert a_{\min ,\lambda ,\beta }-a_{n}\right\Vert _{H^{1}(0,2\alpha
-\varepsilon )}\leq \theta ^{n}\left\Vert w_{\min ,\lambda ,\beta
}-w_{0}\right\Vert _{H^{3}(R)},\quad n=1,2,...,  \label{4.17} \\
& \left\Vert w^{\ast }-w_{n}\right\Vert _{H^{1}(Tr_{\alpha ,\varepsilon
})}\leq C_{1}\xi ^{\rho }+\theta ^{n}\left\Vert w_{\min ,\lambda ,\beta
}-w_{0}\right\Vert _{H^{3}(R)},\quad n=1,2,...,  \label{4.18} \\
& \left\Vert a^{\ast }-a_{n}\right\Vert _{L^{2}(Tr_{\alpha ,\varepsilon
})}\leq C_{1}\xi ^{\rho }+\theta ^{n}\left\Vert w_{\min ,\lambda ,\beta
}-w_{0}\right\Vert _{H^{3}(R)},\quad n=1,2,...  \label{4.19}
\end{align}
\end{theorem}

\begin{remark}
\emph{1. Since the starting point }$W_{0}$\emph{\ of iterations of the
gradient projection method (\ref{4.15}) is an arbitrary point of the ball }$%
B_{0}(2d)$\emph{\ and since the radius }$d>0$\emph{\ of this ball is an
arbitrary number, then estimates (\ref{4.16})-(\ref{4.19}) ensure the global
convergence of the sequence (\ref{4.15}) to the correct solution, see the
first sentence of Introduction.}

\emph{2. We omit below the proofs of Theorem 4.3 and 4.4. Indeed, Theorem
4.3 follows immediately from the combination of Theorem 4.2 with Lemma 2.1
of \cite{bakushinskii2017carleman}. Also, Theorem 4.4 follows immediately
from Theorems 4.2, 4.3, (\ref{4.8}) and (\ref{4.9}). }
\end{remark}

\section{Proofs}

\label{sec:5}

Below in this section $\left( x,t\right) \in R$, where $R$ is the rectangle
defined in (\ref{3.2}).

\subsection{Proof of Theorem 4.1}

\label{sec:5.1}

In this proof $C=C\left( \alpha \right) >0$ denotes different constants
depending only on $\alpha .$ We assume in this proof that the function $u\in
C^{2}\left( \overline{R}\right) \cap H_{0}^{2}\left( R\right) .$ The more
general case $u\in H_{0}^{2}\left( R\right) $ can be obtained from this one
via density arguments. Introduce a new function 
\begin{equation}
v\left( x,t\right) =u\left( x,t\right) e^{-\lambda \left( x+\alpha t\right) }
\label{5.1}
\end{equation}%
and express $u_{xx}-2u_{xt}$ via derivatives of the function $v\left(
x,t\right) .$ We obtain: 
\begin{align*}
& u=ve^{\lambda \left( x+\alpha t\right) },\quad u_{x}=\left( v_{x}+\lambda
v\right) e^{\lambda \left( x+\alpha t\right) },\quad u_{t}=\left(
v_{t}+\lambda \alpha v\right) e^{\lambda \left( x+\alpha t\right) }, \\
& u_{xx}=\left( v_{xx}+2\lambda v_{x}+\lambda ^{2}v\right) e^{\lambda \left(
x+\alpha t\right) },\quad u_{xt}=\left( v_{xt}+\lambda \alpha v_{x}+\lambda
v_{t}+\lambda ^{2}\alpha v\right) e^{\lambda \left( x+\alpha t\right) }, \\
& \left( u_{xx}-2u_{xt}\right) ^{2}e^{-2\lambda \left( x+\alpha t\right) }= 
\left[ \left( v_{xx}-2v_{xt}+\lambda ^{2}\left( 1-2\alpha \right) v\right)
+\left( 2\lambda \left( 1-\alpha \right) v_{x}-2\lambda v_{t}\right) \right]
^{2}.
\end{align*}%
Hence, 
\begin{equation}
\begin{aligned} \left( u_{xx}-2u_{xt}\right) ^{2}e^{-2\lambda \left(
x+\alpha t\right) }&\geq \frac{\left( u_{xx}-2u_{xt}\right) ^{2}e^{-2\lambda
\left( x+\alpha t\right) }}{x+1}\geq \\ &\frac{\left( 4\lambda \left(
1-\alpha \right) v_{x}-4\lambda v_{t}\right) \left( v_{xx}-2v_{xt}+\lambda
^{2}\left( 1-2\alpha \right) v\right) }{x+1}. \label{5.3} \end{aligned}
\end{equation}%
We estimate from below in two steps two products in the second line of (\ref%
{5.3}) involving $v_{x}$ and $v_{t}$.

\textbf{Step 1}. Estimate%
\begin{align*}
&\frac{4\lambda \left( 1-\alpha \right) v_{x}\left( v_{xx}-2v_{xt}+\lambda
^{2}\left( 1-2\alpha \right) v\right) }{x+1} =\left( \frac{2\lambda \left(
1-\alpha \right) v_{x}^{2}}{x+1}\right) _{x}+\frac{2\lambda \left( 1-\alpha
\right) v_{x}^{2}}{\left( x+1\right) ^{2}} + \\
&\left( -\frac{4\lambda \left(1-\alpha \right) v_{x}^{2}}{x+1}\right) _{t}
+\left( \frac{2\lambda ^{3}\left( 1-\alpha \right) \left( 1-2\alpha \right)
v^{2}}{x+1}\right) _{x}+\frac{2\lambda ^{3}\left( 1-\alpha \right) \left(
1-2\alpha \right) v^{2}}{\left( x+1\right) ^{2}}.
\end{align*}

Thus, we have obtained on the first step: 
\begin{equation}
\begin{aligned} &\frac{4\lambda \left( 1-\alpha \right) v_{x}\left(
v_{xx}-2v_{xt}+\lambda ^{2}\left( 1-2\alpha \right) v\right) }{x+1}
=\frac{2\lambda \left( 1-\alpha \right) v_{x}^{2}}{\left( x+1\right)
^{2}}+\frac{2\lambda ^{3}\left( 1-\alpha \right) \left( 1-2\alpha \right)
v^{2}}{\left( x+1\right) ^{2}} + \\ &\left( \frac{2\lambda \left( 1-\alpha
\right) v_{x}^{2}}{x+1}+\frac{2\lambda ^{3}\left( 1-\alpha \right) \left(
1-2\alpha \right) v^{2}}{x+1}\right) _{x}+\left( -\frac{4\lambda \left(
1-\alpha \right) v_{x}^{2}}{x+1}\right) _{t}. \label{5.4} \end{aligned}
\end{equation}

\textbf{Step 2}. Estimate%
\begin{align*}
& -\frac{4\lambda v_{t}\left( v_{xx}-2v_{xt}+\lambda ^{2}\left( 1-2\alpha
\right) v\right) }{x+1}=\left( -\frac{4\lambda v_{t}v_{x}}{x+1}\right) _{x}+%
\frac{4\lambda v_{xt}v_{x}}{x+1}-\frac{4\lambda v_{t}v_{x}}{\left(
x+1\right) ^{2}}+ \\
& \left( \frac{4\lambda v_{t}^{2}}{x+1}\right) _{x}+\frac{4\lambda v_{t}^{2}%
}{\left( x+1\right) ^{2}}+\left( -\frac{2\lambda ^{3}\left( 1-2\alpha
\right) v^{2}}{x+1}\right) _{t}=\frac{4\lambda v_{t}^{2}-4\lambda v_{t}v_{x}%
}{\left( x+1\right) ^{2}}+ \\
& \left( \frac{2\lambda v_{x}^{2}-2\lambda ^{3}\left( 1-2\alpha \right) v^{2}%
}{x+1}\right) _{t}+\left( \frac{4\lambda v_{t}^{2}-4\lambda v_{t}v_{x}}{x+1}%
\right) _{x}.
\end{align*}%
Thus, 
\begin{equation}
\begin{aligned} &-\frac{4\lambda v_{t}\left( v_{xx}-2v_{xt}+\lambda
^{2}\left( 1-2\alpha \right) v\right) }{x+1}=\frac{4\lambda
v_{t}^{2}}{\left( x+1\right) ^{2}}-\frac{4\lambda v_{t}v_{x}}{\left(
x+1\right) ^{2}}\\ &\left( \frac{2\lambda v_{x}^{2}-2\lambda ^{3}\left(
1-2\alpha \right) v^{2}}{x+1}\right) _{t}+\left( \frac{4\lambda
v_{t}^{2}-4\lambda v_{t}v_{x}}{x+1}\right) _{x}. \label{5.5} \end{aligned}
\end{equation}%
Summing up (\ref{5.4}) with (\ref{5.5}) and taking into account (\ref{5.3}),
we obtain%
\begin{equation}
\begin{aligned} \left( u_{xx}-2u_{xt}\right) ^{2}e^{-2\lambda \left(
x+\alpha t\right) }\geq \frac{2\lambda }{\left( x+1\right) ^{2}}\left[
\left( 1-\alpha \right) v_{x}^{2}-2v_{x}v_{t}+2v_{t}^{2}\right] +\\
\frac{2\lambda ^{3}\left( 1-\alpha \right) \left( 1-2\alpha \right)
v^{2}}{\left( x+1\right) ^{2}} +\left( \frac{-2\left( 1-2\alpha \right)
\left( \lambda v_{x}^{2}+\lambda ^{3}v^{2}\right) }{x+1}\right) _{t} \\
\label{5.6} +\left( \frac{2\lambda \left( 1-\alpha \right)
v_{x}^{2}-4\lambda v_{t}v_{x}+4\lambda v_{t}^{2}}{x+1}+\frac{2\lambda
^{3}\left( 1-\alpha \right) \left( 1-2\alpha \right) v^{2}}{x+1}\right)_{x}
\end{aligned}
\end{equation}%
Hence, by Young's inequality 
\begin{equation}
2\lambda \left( 1-\alpha \right) v_{x}^{2}-4\lambda v_{t}v_{x}+4\lambda
v_{t}^{2}\geq 2\lambda \left[ \left( 1-\alpha -\epsilon \right)
v_{x}^{2}+\left( 2-\frac{1}{\epsilon }\right) v_{t}^{2}\right] .
\label{5.60}
\end{equation}

Thus, in order to ensure the positivity of both terms in the right hand side
of (\ref{5.60}), we should have $1/2<\epsilon <1-\alpha .$ We take $\epsilon 
$ as the average of lower and upper bounds of these two inequalities, 
\begin{equation*}
\epsilon =\frac{1}{2}\left( \frac{1}{2}+\left( 1-\alpha \right) \right) =%
\frac{3-2\alpha }{4}.
\end{equation*}%
Hence, (\ref{5.60}) becomes 
\begin{equation}
2\lambda \left( 1-\alpha \right) v_{x}^{2}-4\lambda v_{t}v_{x}+4\lambda
v_{t}^{2}\geq \frac{\lambda \left( 1-2\alpha \right) }{2}v_{x}^{2}+\frac{%
4\lambda \left( 1-2\alpha \right) }{3-2\alpha }v_{t}^{2}.  \label{5.7}
\end{equation}

Note that since $u\in C^{2}\left( \overline{R}\right) \cap H_{0}^{2}\left(
R\right) ,$ then by (\ref{5.1}) $v\left( 0,t\right) =v_{x}\left( 0,t\right)
=0.$ Hence, integrating (\ref{5.6}) over $R$ and taking into account (\ref%
{5.7}), we obtain 
\begin{equation}
\begin{aligned} &\int\displaylimits_{R}\left( u_{xx}-2u_{xt}\right)
^{2}e^{-2\lambda \left( x+\alpha t\right) } \geq C\lambda
\int\displaylimits_{R}\left( v_{x}^{2}+v_{t}^{2}\right) dxdt+C\lambda
^{3}\int\displaylimits_{R}v^{2}dxdt \\ &+C\lambda
\int\displaylimits_{0}^{1}v_{x}^{2}\left( x,0\right) dx +
C\lambda^{3}\int\displaylimits_{0}^{1}v^{2}\left( x,0\right) dx -C\lambda
\int\displaylimits_{0}^{1}v_{x}^{2}\left( x,T\right) dx-C\lambda
^{3}\int\displaylimits_{0}^{1}v^{2}\left( x,T\right) dx. \label{5.8}
\end{aligned}
\end{equation}%
We now replace in (\ref{5.8}) the function $v$ with the function $u$ via (%
\ref{5.1}). We have 
\begin{align*}
& \lambda v_{x}^{2}=\lambda \left( u_{x}^{2}-2\lambda u_{x}u+\lambda
^{2}u^{2}\right) e^{-2\lambda \left( x+\alpha t\right) }\geq \left( \frac{%
\lambda }{2}u_{x}^{2}-\lambda ^{3}u^{2}\right) e^{-2\lambda \left( x+\alpha
t\right) }, \\
& \lambda v_{t}^{2}=\lambda \left( u_{t}^{2}-2\lambda \alpha u_{t}u+\lambda
^{2}\alpha ^{2}u^{2}\right) e^{-2\lambda \left( x+\alpha t\right) }\geq
\left( \frac{\lambda }{2}u_{t}^{2}-\lambda ^{3}\alpha ^{2}u^{2}\right)
e^{-2\lambda \left( x+\alpha t\right) }.
\end{align*}%
Thus, 
\begin{equation*}
C\lambda \left( v_{x}^{2}+v_{t}^{2}\right) \geq \frac{C}{4}\lambda \left(
v_{x}^{2}+v_{t}^{2}\right) \geq \left( \frac{C}{8}\lambda \left(
u_{x}^{2}+u_{t}^{2}\right) -\frac{C}{2}\lambda ^{3}u^{2}\right) e^{-2\lambda
\left( x+\alpha t\right) }.
\end{equation*}

\noindent Hence, (\ref{5.8}) implies the following estimate, which is
equivalent with (\ref{4.1}): 
\begin{align*}
&\int \displaylimits_{R}\left( u_{xx}-2u_{xt}\right) ^{2}e^{-2\lambda \left(
x+\alpha t\right) }\geq \frac{C}{8}\lambda \int \displaylimits_{R}\left(
u_{x}^{2}+u_{t}^{2}\right) e^{-2\lambda \left( x+\alpha t\right) }dxdt \\
&+\frac{C}{2}\lambda ^{3}\int \displaylimits_{R}u^{2}e^{-2\lambda \left(
x+\alpha t\right) }dxdt+\frac{C}{8}\lambda \int \displaylimits%
_{0}^{1}u_{x}^{2}\left( x,0\right) e^{-2\lambda x}dx \\
&+ \frac{C}{2}\lambda ^{3}\int \displaylimits_{0}^{1}u^{2}\left( x,0\right)
e^{-2\lambda x}dx-C\lambda e^{-2\lambda \alpha T}\int \displaylimits%
_{0}^{1}u_{x}^{2}\left( x,T\right) dx-C\lambda ^{3}e^{-2\lambda \alpha
T}\int \displaylimits_{0}^{1}u^{2}\left( x,T\right) dx.\text{ \ }\square
\end{align*}

\subsection{Proof of Theorem 4.2}

\label{sec:5.2}

Let two arbitrary functions $w_{1},w_{2}\in \overline{B\left(
d,p_{0},p_{1}\right) }$. Denote $h=w_{2}-w_{1}.$ Then $h\in \overline{%
B_{0}\left( 2d\right) }.\label{5.110}$ Note that embedding theorem implies
that sets $\overline{B\left( d,p_{0},p_{1}\right) },\overline{B_{0}\left(
2d\right) }\subset C^{1}\left( \overline{R}\right) $, 
\begin{equation}
\left\Vert w\right\Vert _{C^{1}\left( \overline{R}\right) }\leq C_{1},\quad
\forall w\in \overline{B\left( d,p_{0},p_{1}\right) },\quad \left\Vert
h\right\Vert _{C^{1}\left( \overline{R}\right) }\leq C_{1}.  \label{5.11}
\end{equation}%
It follows from (\ref{3.4}) that in this proof, we should first estimate
from below $\left[ L\left( w_{1}+h\right) \right] ^{2}-\left[ L\left(
w_{1}\right) \right] ^{2}.$ We will single out the linear and nonlinear
parts, with respect to $h$, of this expression. By (\ref{2.160}): 
\begin{equation}
\begin{aligned} &L\left( w_{1}+h\right) =L\left( w_{1}\right)
+h_{xx}-2h_{xt}+2h_{x}\int\displaylimits_{0}^{t}w_{1x}\left( x,\tau \right)
d\tau +2w_{1x}\int\displaylimits_{0}^{t}h_{x}\left( x,\tau \right) d\tau \\
&-2h_{x}w_{1} -2h_{x}h -2w_{1x}h
-2h_{t}\int\displaylimits_{0}^{t}w_{1x}\left( x,\tau \right) d\tau
-2w_{1t}\int\displaylimits_{0}^{t}h_{x}\left( x,\tau \right) d\tau \\
&+2\left[ h_{x}\int\displaylimits_{0}^{t}h_{x}\left( x,\tau \right) d\tau
-h_{t}\int\displaylimits_{0}^{t}h_{x}\left( x,\tau \right) d\tau \right] =
L\left( w_{1}\right) +L_{lin}\left( h\right) +L_{nl}\left( h\right) ,
\label{5.12} \end{aligned}
\end{equation}%
where $L_{lin}\left( h\right) $ and $L_{nl}\left( h\right) $ are linear and
nonlinear, with respect to $h$, parts of (\ref{5.12}), and their forms are
clear from (\ref{5.12}). Hence,%
\begin{equation}
\begin{aligned} \left[ L\left( w_{1}+h\right) \right] ^{2}&-\left[ L\left(
w_{1}\right) \right] ^{2}=2L\left( w_{1}\right) L_{lin}\left( h\right)
+\left( L_{lin}\left( h\right) \right) ^{2}+\\ &\left( L_{nl}\left( h\right)
\right) ^{2}+2L_{lin}\left( h\right) L_{nl}\left( h\right) +2L\left(
w_{1}\right) L_{nl}\left( h\right) . \label{5.13} \end{aligned}
\end{equation}

Using (\ref{5.11}), (\ref{5.12}) and the Cauchy-Schwarz inequality, we obtain%
\begin{eqnarray}
&&\left( L_{lin}\left( h\right) \right) ^{2}+\left( L_{nl}\left( h\right)
\right) ^{2}+2L_{lin}\left( h\right) L_{nl}\left( h\right) +2L\left(
w_{1}\right) L_{nl}\left( h\right)   \label{5.14} \\
&\geq &\frac{1}{2}\left( h_{xx}-2h_{xt}\right) ^{2}-C_{1}\left[
h_{x}^{2}+h_{t}^{2}+h^{2}+\left( \int \displaylimits_{0}^{t}h_{x}\left(
x,\tau \right) d\tau \right) ^{2}\right] .  \notag
\end{eqnarray}%
Let $\left( \cdot ,\cdot \right) $ denotes the scalar product in $%
H^{3}\left( R\right) .$ It follows from (\ref{3.4}) and (\ref{5.13}) that 
\begin{equation}
J_{\lambda ,\beta }\left( w_{1}+h\right) -J_{\lambda ,\beta }\left(
w_{1}\right) =A\left( h\right) +B\left( h\right) ,  \label{5.15}
\end{equation}%
where $A\left( h\right) :H_{0}^{3}\left( R\right) \rightarrow \mathbb{R}$ is
a bounded linear functional, 
\begin{equation*}
A\left( h\right) =\int \displaylimits_{R}2L\left( w_{1}\right) L_{lin}\left(
h\right) \varphi _{\lambda }dxdt+2\beta \left( w_{1},h\right) 
\end{equation*}%
and $B\left( h\right) $ is a nonlinear functional,%
\begin{equation}
B\left( h\right) =  \label{5.150}
\end{equation}%
\begin{equation*}
\int\displaylimits_{R}\left[ \left( L_{lin}\left( h\right) \right) ^{2}+\left(
L_{nl}\left( h\right) \right) ^{2}+2L_{lin}\left( h\right) L_{nl}\left(
h\right) +2L\left( w_{1}\right) L_{nl}\left( h\right) \right] \varphi
_{\lambda }dxdt+\beta \left\Vert h\right\Vert _{H^{3}\left( R\right) }^{2}.
\end{equation*}%
By the Riesz theorem, there exists unique point $\widetilde{A}\in
H_{0}^{3}\left( R\right) $ such that 
\begin{equation}
A\left( h\right) =\left( \widetilde{A},h\right) ,\hspace{0.3em}\forall h\in
H_{0}^{3}\left( R\right) .  \label{6}
\end{equation}%
Next, it follows from (\ref{5.15})-(\ref{6}) that 
\begin{equation*}
\lim_{\left\Vert h\right\Vert _{H^{3}\left( R\right) }\rightarrow 0}\frac{%
\left\vert J_{\lambda ,\beta }\left( w_{1}+h\right) -J_{\lambda ,\beta
}\left( w_{1}\right) -\left( \widetilde{A},h\right) \right\vert }{\left\Vert
h\right\Vert _{H^{3}\left( R\right) }}=0.
\end{equation*}

Hence, $\widetilde{A}\in H_{0}^{3}\left( R\right) $ is the Fr\'{e}chet
derivative $J_{\lambda ,\beta }^{\prime }\left( w_{1}\right) \in
H_{0}^{3}\left( R\right) $ of the functional $J_{\lambda ,\beta }\left(
w_{1}\right) $ at the point $w_{1},$%
\begin{equation}
\widetilde{A}=J_{\lambda ,\beta }^{\prime }\left( w_{1}\right) .  \label{7}
\end{equation}%
Next, (\ref{3.4}) and (\ref{5.14})-(\ref{7}) imply that for all $\lambda
\geq 1$%
\begin{equation}
\begin{aligned} &J_{\lambda ,\beta }\left( w_{1}+h\right) -J_{\lambda ,\beta
}\left( w_{1}\right) -J_{\lambda ,\beta }^{\prime }\left( w_{1}\right)
\left( h\right) \geq \frac{1}{2}\int\displaylimits_{R}\left(
h_{xx}-2h_{xt}\right) ^{2}\varphi _{\lambda }dxdt \label{5.16} \\
&-C_{1}\int\displaylimits_{R}\left[ h_{x}^{2}+h_{t}^{2}+h^{2}+\left(
\int\displaylimits_{0}^{t}h_{x}\left( x,\tau \right) d\tau \right)
^{2}\right] \varphi _{\lambda }dxdt+\beta \left\Vert h\right\Vert
_{H^{3}\left( R\right) }^{2}. \end{aligned}
\end{equation}

Combining Lemma 3.1, Theorem 4.1 and (\ref{5.16}) and also assuming that $%
\lambda \geq \lambda _{0},$ we obtain 
\begin{equation}
\begin{aligned} &J_{\lambda ,\beta }\left( w_{1}+h\right) -J_{\lambda ,\beta
}\left( w_{1}\right) -J_{\lambda ,\beta }^{\prime }\left( w_{1}\right)
\left( h\right) \geq C\lambda \int\displaylimits_{R}\left(
h_{x}^{2}+h_{t}^{2}\right) \varphi _{\lambda }dxdt \\
&+C\lambda^{3}\int\displaylimits_{R}h^{2}\varphi _{\lambda }dxdt +\beta
\left\Vert h\right\Vert _{H^{3}\left( R\right)
}^{2}-C_{1}\int\displaylimits_{R}\left( h_{x}^{2}+h_{t}^{2}+h^{2}\right)
\varphi _{\lambda }dxdt \\ &+C\lambda
\int\displaylimits_{0}^{1}h_{x}^{2}\left( x,0\right) e^{-2\lambda x}dx
+C\lambda ^{3}\int\displaylimits_{0}^{1}h^{2}\left( x,0\right) e^{-2\lambda
x}dx \\ &-C\lambda e^{-2\lambda \alpha
T}\int\displaylimits_{0}^{1}h_{x}^{2}\left( x,T\right) dx-C\lambda
^{3}e^{-2\lambda \alpha T}\int\displaylimits_{0}^{1}h^{2}\left( x,T\right)
dx. \label{5.17} \end{aligned}
\end{equation}%
Choose $\lambda _{1}=\lambda _{1}\left( \alpha ,\varepsilon ,d\right) \geq $%
\emph{\ }$\lambda _{0}\geq 1$ so large that $C\lambda _{1}>2C_{1}$ and then
take in (\ref{5.17}) $\lambda \geq \lambda _{1}.$ We obtain%
\begin{equation}
\begin{aligned} &J_{\lambda ,\beta }\left( w_{1}+h\right) -J_{\lambda ,\beta
}\left( w_{1}\right) -J_{\lambda ,\beta }^{\prime }\left( w_{1}\right)
\left( h\right) \geq C_{1}\lambda \int\displaylimits_{R}\left(
h_{x}^{2}+h_{t}^{2}\right) \varphi _{\lambda }dxdt\\ &+C_{1}\lambda
^{3}\int\displaylimits_{R}h^{2}\varphi _{\lambda }dxdt + C_{1}\lambda
\int\displaylimits_{0}^{1}h_{x}^{2}\left( x,0\right) e^{-2\lambda
x}dx+C_{1}\lambda ^{3}\int\displaylimits_{0}^{1}h^{2}\left( x,0\right)
e^{-2\lambda x}dx \\ &+\beta \left\Vert h\right\Vert _{H^{3}\left( R\right)
}^{2} - C_{1}\lambda e^{-2\lambda \alpha
T}\int\displaylimits_{0}^{1}h_{x}^{2}\left( x,T\right) dx-C_{1}\lambda
^{3}e^{-2\lambda \alpha T}\int\displaylimits_{0}^{1}h^{2}\left( x,T\right)
dx. \label{5.18} \end{aligned}
\end{equation}%
Since $Tr_{\alpha ,\varepsilon }\subset Tr\subset R$ and since the interval $%
\left( 0,2\alpha -\varepsilon \right) \subset \left( 0,1\right) $ and also
since $\varphi _{\lambda }\left( x,t\right) \geq e^{-2\lambda \left( 2\alpha
-\varepsilon \right) }$ in $Tr_{\alpha ,\varepsilon },$ then we obtain from (%
\ref{5.18})%
\begin{equation*}
\begin{aligned} &J_{\lambda ,\beta }\left( w_{1}+h\right) -J_{\lambda ,\beta
}\left( w_{1}\right) -J_{\lambda ,\beta }^{\prime }\left( w_{1}\right)
\left( h\right) \geq C_{1}e^{-2\lambda \left( 2\alpha -\varepsilon \right)
}\left\Vert h\right\Vert _{H^{1}\left( Tr_{\alpha ,\varepsilon }\right)
}^{2}+ \\ &C_{1}e^{-2\lambda \left( 2\alpha -\varepsilon \right) }\left\Vert
h\left( x,0\right) \right\Vert _{H^{1}\left( 0,2\alpha -\varepsilon \right)
}^{2} +\beta \left\Vert h\right\Vert _{H^{3}\left( R\right)
}^{2}-C_{1}\lambda ^{3}e^{-2\lambda \alpha T}\left\Vert h\left( x,T\right)
\right\Vert _{H^{1}\left( Tr_{\alpha ,\varepsilon }\right) }^{2},
\hspace{0.3em} \forall \lambda \geq \lambda _{1}. \end{aligned}
\end{equation*}%
By the trace theorem $\left\Vert h\left( x,T\right) \right\Vert
_{H^{1}\left( 0,2\alpha -\varepsilon \right) }^{2}\leq C_{1}\left\Vert
h\right\Vert _{H^{3}\left( R\right) }^{2}.$ Hence, taking $\beta \in \left[
2e^{-\lambda \alpha T},1\right) ,$ we obtain the following estimate for all $%
\lambda \geq \lambda _{1}$: 
\begin{equation}
\begin{aligned} J_{\lambda ,\beta }\left( w_{1}+h\right) -J_{\lambda ,\beta
}\left( w_{1}\right) -J_{\lambda ,\beta }^{\prime }\left( w_{1}\right)
\left( h\right) \geq C_{1}e^{-2\lambda \left( 2\alpha -\varepsilon \right)
}\left\Vert h\right\Vert _{H^{1}\left( Tr_{\alpha ,\varepsilon }\right)
}^{2} \\ +C_{1}e^{-2\lambda \left( 2\alpha -\varepsilon \right) }\left\Vert
h\left( x,0\right) \right\Vert _{H^{1}\left( 0,2\alpha -\varepsilon \right)
}^{2} +\frac{\beta }{2}\left\Vert h\right\Vert _{H^{3}\left( R\right)
}^{2}.\quad \end{aligned}
\end{equation}

This estimate is equivalent with our target estimate (\ref{4.3}). \ \ $%
\square $

\subsection{Proof of Theorem 4.5}

\label{sec:5.3}

Let $\lambda \geq \lambda _{2}.$ Temporary denote $I_{\lambda ,\beta }\left(
W,G\right) :=J_{\lambda ,\beta }\left( W+G\right) .$ Consider $I_{\lambda
,\beta }\left( W^{\ast },G\right) ,$ 
\begin{equation}
\begin{aligned} I_{\lambda ,\beta }\left( W^{\ast },G\right) =J_{\lambda
,\beta }\left( W^{\ast }+G\right) =\int \displaylimits_{R}\left[ L\left(
W^{\ast }+G\right) \right] ^{2}\varphi _{\lambda }dxdt+\beta \left\Vert
W^{\ast }+G\right\Vert _{H^{3}\left( R\right) }^{2} = \label{5.19} \\
J_{\lambda ,\beta }^{0}\left( W^{\ast }+G\right) +\beta \left\Vert W^{\ast
}+G\right\Vert _{H^{3}\left( R\right) }^{2} \end{aligned}
\end{equation}
Since $L\left( W^{\ast }+G^{\ast }\right) =L\left( w^{\ast }\right) =0,$
then 
\begin{equation*}
L\left( W^{\ast }+G\right) =L\left( W^{\ast }+G^{\ast }+\left( G-G^{\ast
}\right) \right) =L\left( W^{\ast }+G^{\ast }\right) +\widehat{L}\left(
G-G^{\ast }\right) =\widehat{L}\left( G-G^{\ast }\right) ,
\end{equation*}%
where by (\ref{2.160}) and (\ref{4.5}), $\left\vert \widehat{L}\left(
G-G^{\ast }\right) \left( x,t\right) \right\vert \leq C_{1}\xi $ for all $%
\left( x,t\right) \in \overline{R}.$ Hence, by (\ref{5.19}) 
\begin{equation}
I_{\lambda ,\beta }\left( W^{\ast },G\right) \leq C_{1}\left( \xi ^{2}+\beta
\right) .  \label{5.20}
\end{equation}%
We have 
\begin{equation}
W^{\ast }-W_{\min ,\lambda ,\beta }=\left( W^{\ast }+G\right) -\left(
W_{\min ,\lambda ,\beta }+G\right) =\left( w^{\ast }-w_{\min ,\lambda ,\beta
}\right) +\left( G-G^{\ast }\right) .  \label{5.200}
\end{equation}

Also, by (\ref{4.5}) and the trace theorem 
\begin{equation}
\left\Vert G\left( x,0\right) -G^{\ast }\left( x,0\right) \right\Vert
_{H^{1}\left( 0,2\alpha -\varepsilon \right) }\leq C_{1}\xi .  \label{5.201}
\end{equation}%
Hence, (\ref{4.5}), (\ref{5.200}) and (\ref{5.201}) imply 
\begin{equation*}
\begin{aligned} &\left\Vert W^{\ast }-W_{\min ,\lambda ,\beta }\right\Vert
_{H^{1}\left( Tr_{\alpha ,\varepsilon }\right) }^{2}\geq
\frac{1}{2}\left\Vert w^{\ast }-w_{\min ,\lambda ,\beta }\right\Vert
_{H^{1}\left( Tr_{\alpha ,\varepsilon }\right) }^{2}-C_{1} \xi^{2}, \\
&\left\Vert W^{\ast }\left( x,0\right) -W_{\min ,\lambda ,\beta }\left(
x,0\right) \right\Vert _{H^{1}\left( 0,2\alpha -\varepsilon \right)
}^{2}\geq \frac{1}{2}\left\Vert w^{\ast }-w_{\min ,\lambda ,\beta
}\right\Vert _{H^{1}\left( 0,2\alpha -\varepsilon \right) }^{2}-C_{1} \xi^{2} \\
\frac{\beta }{2}&\left\Vert W^{\ast }-W_{\min ,\lambda
,\beta }\right\Vert _{H^{3}\left( R\right) }^{2}\geq \frac{\beta
}{4}\left\Vert w^{\ast }-w_{\min ,\lambda ,\beta }\right\Vert _{H^{3}\left(
R\right) }^{2}-\frac{\beta }{2}\xi ^{2} \end{aligned}
\end{equation*}%
Hence, using (\ref{4.10}), we obtain 
\begin{equation}
\begin{aligned} &I_{\lambda ,\beta }\left( W^{\ast },G\right) -I_{\lambda
,\beta }\left( W_{\min ,\lambda ,\beta },G\right) -I_{\lambda ,\beta
}^{\prime }\left( W_{\min ,\lambda ,\beta },G\right) \left( W^{\ast
}-W_{\min ,\lambda ,\beta }\right) \geq \\ &C_{1}e^{-2\lambda \left( 2\alpha
-\varepsilon \right) }\left\Vert w^{\ast }-w_{\min ,\lambda ,\beta
}\right\Vert _{H^{1}\left( Tr_{\alpha ,\varepsilon }\right) }^{2}-C_{1} \xi^{2}
\label{5.21} \\ +&C_{1}e^{-2\lambda \left( 2\alpha -\varepsilon \right)
}\left\Vert w^{\ast }\left( x,0\right) -w_{\min ,\lambda ,\beta }\left(
x,0\right) \right\Vert _{H^{1}\left( 0,2\alpha -\varepsilon \right)
}^{2}. \end{aligned}
\end{equation}%
By (\ref{4.11}) 
\begin{equation*}
-I_{\lambda ,\beta }^{\prime }\left( W_{\min ,\lambda ,\beta },G\right)
\left( W^{\ast }-W_{\min ,\lambda ,\beta }\right) \leq 0.
\end{equation*}
Hence,%
\begin{equation*}
I_{\lambda ,\beta }\left( W^{\ast },G\right) -I_{\lambda ,\beta }\left(
W_{\min ,\lambda ,\beta },G\right) -I_{\lambda ,\beta }^{\prime }\left(
W_{\min ,\lambda ,\beta },G\right) \left( W^{\ast }-W_{\min ,\lambda ,\beta
}\right) \leq I_{\lambda ,\beta }\left( W^{\ast },G\right) .
\end{equation*}%
Comparing this with (\ref{5.20}) with (\ref{5.21}) and dropping the term
with $\beta $ in (\ref{5.21}), we obtain%
\begin{equation}
e^{-2\lambda \left( 2\alpha -\varepsilon \right) }\left( \left\Vert w^{\ast
}-w_{\min ,\lambda ,\beta }\right\Vert _{H^{1}\left( Tr_{\alpha ,\varepsilon
}\right) }^{2}+\left\Vert w^{\ast }\left( x,0\right) -w_{\min ,\lambda
,\beta }\left( x,0\right) \right\Vert _{H^{1}\left( 0,2\alpha -\varepsilon
\right) }^{2}\right)   \label{200}
\end{equation}%
\begin{equation*}
\leq C_{1}\left( \xi ^{2}+\beta \right) .
\end{equation*}%
Dividing both sides of (\ref{200}) by $e^{-2\lambda \left( 2\alpha
-\varepsilon \right) }$ and recalling that by (\ref{4.121}) $\beta
=2e^{-\lambda \alpha T}$, we obtain%
\begin{equation*}
\left\Vert w^{\ast }-w_{\min ,\lambda ,\beta }\right\Vert _{H^{1}\left(
Tr_{\alpha ,\varepsilon }\right) }^{2}+\left\Vert w^{\ast }\left( x,0\right)
-w_{\min ,\lambda ,\beta }\left( x,0\right) \right\Vert _{H^{1}\left(
0,2\alpha -\varepsilon \right) }^{2}
\end{equation*}%
\begin{equation}
\leq C_{1}\xi ^{2}e^{2\lambda \left( 2\alpha -\varepsilon \right)
}+C_{1}\exp \left( -\lambda \left( \alpha \left( T-4\right) +2\varepsilon
\right) \right) .  \label{5.230}
\end{equation}

Since $T\geq 4,$ then $-\lambda \left( \alpha \left( T-4\right)
+2\varepsilon \right) <0.$ Since we have chosen $\lambda =\lambda \left( \xi
\right) $ and $\beta =\beta \left( \xi \right) $ as in (\ref{4.121}), then
in (\ref{5.230}) $\xi ^{2}e^{2\lambda \left( 2\alpha -\varepsilon \right)
}=\xi $ and $\exp \left( -\lambda \left( \alpha \left( T-4\right)
+2\varepsilon \right) \right) =\xi ^{\sigma }.$ Hence, target estimates (\ref%
{4.14}) follow from (\ref{2.19}), (\ref{4.12}) and (\ref{5.230}). $\ \ \
\square $

\subsection{Proof of Theorem 4.6}

\label{sec:5.4}

The existence of the number $\theta \in \left( 0,1\right) $ as well as
convergence rates (\ref{4.16}) and (\ref{4.17}) follow immediately from a
combination of Theorem 4.2 with Theorem 2.1 of \cite%
{bakushinskii2017carleman}. Convergence rate (\ref{4.18}) follows
immediately from the triangle inequality, (\ref{4.14}) and (\ref{4.16}).
Similarly, convergence rate (\ref{4.19}) follows immediately from the
triangle inequality, (\ref{4.14}) and (\ref{4.17}). \ $\square $

\section{Numerical Implementation}

\label{sec:6}

To computationally simulate the data (\ref{1.4}) for our CIP, we solve the
forward problem (\ref{2.180})-(\ref{2.182}) by the finite difference method
in the domain $\{(x,t)\in (-A,A)\times (0,T)\}$. In all our computations of
the forward problem (\ref{2.180})-(\ref{2.182}) we take $A=B=2.2,T=4$. For a
given function $a(x)$ we compute the solution $u_{i,j}=u(x_{i},t_{j})$ on
the rectangular mesh with $N_{x}=1024$ spatial and $N_{t}=1024$ temporal
grid points.

Now, even though Theorems 4.5 and 4.6 work only for $T\geq 4,$ we use $T=2$
in our computations of the inverse problem. Also, when computing the inverse
problem, we take $A=1.1.$ Thus, the rectangle $R$ in (\ref{3.2}) is replaced
in our computations of the inverse problem with the rectangle $R^{\prime }$, 
\begin{equation*}
R^{\prime }=\left( 0,A\right) \times \left( 0,T\right) =(0,1.1)\times (0,2).
\end{equation*}%
In order to avoid the inverse crime, we work in the inverse problem with the
rectangular mesh of $N_{x}\times N_{t}=60\times 50$ grid points. The
absorbing boundary condition (\ref{2.182}) at $x=A$ gives us the following
direct analog of boundary condition (\ref{6.4}): 
\begin{equation}
w_{x}\left( 1.1,t\right) =0.  \label{8}
\end{equation}%
We have observed numerically that this condition provides a better stability
for our computations of the inverse problem, as compared with the case when
condition (\ref{8}) is absent.

The finite difference approximations of differential operators in (\ref%
{2.160}) are used on the rectangular mesh with $h=(h_{x},h_{t})$. Denote $%
w(x_{i},t_{j})=w^{i,j}$. We write \ the functional $J_{\lambda ,\beta
}\left( w\right) $ in (\ref{3.4}) in the finite difference form as: 
\begin{equation}
\begin{aligned} J_{\lambda,\beta, \mu}^{h} (w^{i,j}) &=
h_{x}h_{t}\sum_{i=3}^{N_{x}-1}\sum_{j=1}^{N_{t}-1}\Bigg(\frac{w^{i,j}-2
w^{i+1,j} +w^{i+2,j}}{h_x^2}\\
&-2\frac{w^{i+1,j+1}-w^{i+1,j}-w^{i,j+1}+w^{i,j}}{h_x h_t} \\&+2h_t
\frac{w^{i+1,j}-w^{i,j}}{h_x}\sum_{l=1}^{N_{t}-1}
\left(\frac{w^{i+1,l}-w^{i,l}}{h_x} \right)
-2\frac{w^{i+1,j}-w^{i,j}}{h_x}w^{i,j}
\\&-2(w^{i,j+1}-w^{i,j})\sum_{l=1}^{N_{t}-1}
\left(\frac{w^{i+1,l}-w^{i,l}}{h_x} \right) \Bigg)^{2} e^{ -2\lambda (
x_i+\alpha t_j)}\\ &+\beta h_{x} h_{t}
\sum_{i=3}^{N_{x}-1}\sum_{j=1}^{N_{t}-1}\Bigg(\left( w^{i,j} \right)^2
+\left( \frac{w^{i+1,j}-w^{i,j}}{h_x} \right)^2
+\left(\frac{w^{i,j+1}-w^{i,j}}{h_t} \right)^2 \\ &+ \left( \frac{w^{i,j}-2
w^{i+1,j} +w^{i+2,j}}{h_x^2} \right)^2 +\left( \frac{w^{i,j}-2 w^{i,j+1}
+w^{i,j+2}}{h_t^2} \right)^2 \Bigg) \\ &+\mu \sum_{j=1}^{N_{t}-1} \left(
\frac{w^{Nx,j} - w^{Nx-1,j}}{h_x} \right)^2. \end{aligned}  \label{6.5}
\end{equation}%
Next, we minimize functional (\ref{6.5}) with respect to the values $w^{i,j}$
of the unknown function $w\left( x,t\right) $ at grid points $(x_{i},t_{j})$%
. To speed up computations, the gradient of the functional (\ref{6.5}) is
written in an explicit form, using Kronecker symbols, as in \cite%
{klibanov2008new}. For brevity, we do not bring in these formulas here.

\begin{remark}
\emph{1. In fact the functional (\ref{6.5}), which is used to conduct
numerical studies, is a slightly modified finite difference version of (\ref%
{3.4}). In our computations, we took the Tikhonov regularization term in the
finite difference analog of $H^{2}\left( R^{\prime }\right) $ instead of $%
H^{3}\left( R^{\prime }\right) $. Note that since the number of grid points
is not exceedingly large here ($N_{x}=60,N_{t}=50$), then all discrete norms
are basically equivalent. Additionally, the boundary term with the
coefficient }$\mu >>1$\emph{\ is added in (\ref{6.5}) to ensure that the
minimizer satisfies boundary condition (\ref{8}).}

\emph{\ 2. We choose parameters $\lambda ,\alpha ,\beta $ and }$\mu $\emph{\
so that the numerical method provides a good reconstruction of a reference
function $a(x)$ of our choice depicted on Figure 3(A). The values of our
parameters were found by the trial and error procedure. It is important
though that exactly the same values of those parameters were used then in
three subsequent tests. Those values were: 
\begin{equation}
\lambda =2,\hspace{0.3em}\alpha =1/2,\hspace{0.3em}\beta =10^{-4},\hspace{%
0.3em}\mu =10^{2}.  \label{9}
\end{equation}%
We note that even though the parameter $\lambda $ has to be sufficiently
large, $\lambda =2$ worked quite well in our numerical experiments. This is
similar with all above cited works about numerical studies of the
convexification. The topic of optimal choices of these parameters is outside
of the scope of this paper. Also, see below a brief discussion of the choice
of parameters }$\lambda $ and $\alpha .$

\emph{3. Even though Theorem 4.6 guarantees the global convergence of the
gradient projection method, we have observed in our computations that just
the straightforward gradient descent method works well. This method is
simpler to implement than the gradient projection method since one does not
need to use the orthogonal projection operator }$P_{B_{0}}$\emph{\ in (\ref%
{4.15}). Thus, we have not subtracted the function }$G$\emph{\ from the
function }$w$\emph{\ and minimized, therefore, the functional }$J_{\lambda
,\beta }$\emph{\ instead of the functional }$I_{\lambda ,\beta }.$\emph{\ In
other words, (\ref{4.15}) was replaced with}%
\begin{equation}
w_{n}=w_{n-1}-\gamma J_{\lambda ,\beta }^{\prime }(w_{n-1})),\quad
n=1,2,\dots   \label{6.50}
\end{equation}%
\emph{Note that }$J_{\lambda ,\beta }^{\prime }\in H_{0}^{3}\left( R^{\prime
}\right) .$\emph{\ This means that all functions }$w_{n}$\emph{\ of the
sequence (\ref{6.50}) satisfy the same boundary conditions }$p_{0},p_{1}$ 
\emph{(\ref{2.17})}.\emph{\ \ We took }$\gamma =10^{-5}$ \emph{at the first
step of the gradient descent method and adjusted it using line search at
every subsequent iteration.}

\emph{4. We choose the starting point $w_{0}(x,t)$ of the process (\ref{6.50}%
) as $w_{0}(x,t)=-(p_{1}(t)x^{2})/2.2+p_{1}(t)x+p_{0}(t)$.} \emph{It is easy
to see that the function $w_{0}(x,t)$ satisfies boundary conditions (\ref%
{2.17}) as well as boundary condition (\ref{8}). Hence, we set at the first
step of the minimization procedure} 
\begin{equation*}
a_{0}(x)=2(w_{0})_{x}(x,0)=2p_{1}(0)(1-2x/2.2).
\end{equation*}%
\emph{In most cases }$p_{1}(0)=0$\emph{, which means that the initial
function }$a_{0}(x)\equiv 0$ in most cases\emph{. Using (\ref{2.19}), we set 
$a_{n}(x)=2(w_{n})_{x}(x,0)$, where the function $w_{n}(x,t)$ is computed on
the $n$-th step of the minimization procedure.}

\emph{5. The stopping criterion for our minimization process is 
\begin{equation*}
\Vert a_{n+1}-a_{n}\Vert _{L^{2}\left( 0,1\right) }/\Vert a_{n}\Vert
_{L^{2}\left( 0,1\right) }\hspace{0.3em}\leq \hspace{0.3em}10^{-2}.
\end{equation*}%
}
\end{remark}

\begin{figure}[tbp]
\begin{center}
\subfloat[\emph{$\Vert J^h_{\lambda,\beta, \mu}(w_n) \Vert_{\infty}$ for
$n = 1, \dots 30$.}]{\includegraphics[width =.48\textwidth]{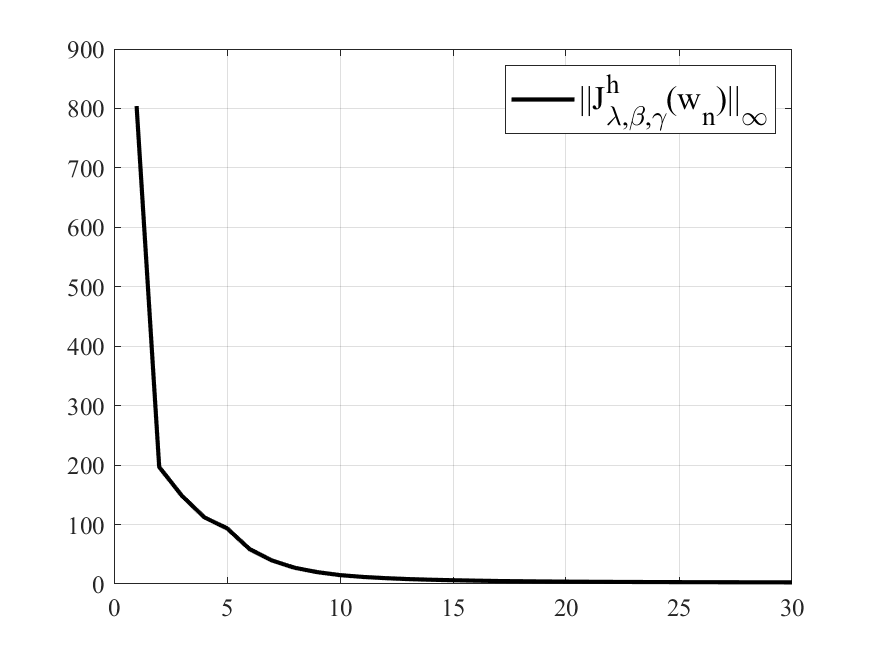}} 
\subfloat[\emph{$u(0,t)$ and $u^{\xi}(0,t)$}, $\xi =
0.1$.]{\includegraphics[width =.48\textwidth]{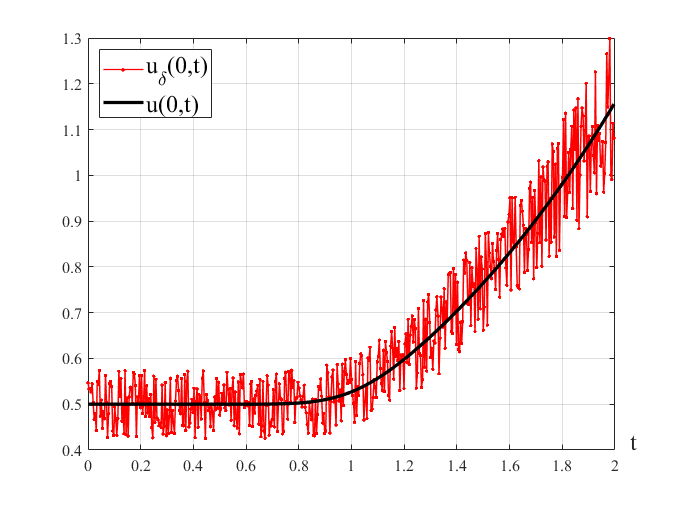}} \quad 
\subfloat[\emph{$u_x(0,t)$ and $u_x^{\xi}(0,t)$}, $\xi =
0.1$.]{\includegraphics[width =.48\textwidth]{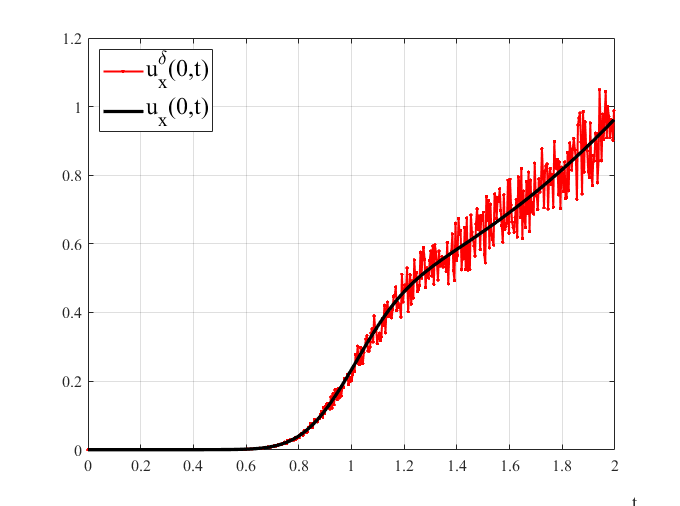}} \quad %
\subfloat[\emph{$\left\Vert
J_{0,\beta ,\mu }\left( w_{n}\right) \right\Vert _{\infty }$ for $n=1,...,10$. }]{\includegraphics[width =.48\textwidth]{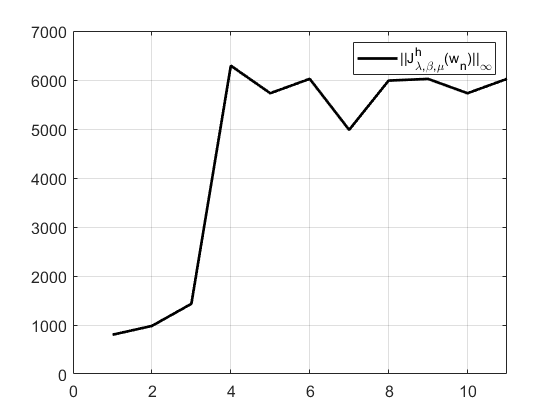}}
\end{center}
\caption{\emph{The comparison of noiseless and noisy data. Figure 2(A) shows
the norm of the functional (\protect\ref{6.5}) for each iteration of the
gradient descent for the test function depicted on Figure 3(A). Figure 2(D) corresponds to our test for $\lambda$ = 0; see the text.}}
\label{fig2}
\end{figure}

\subsection{Data pre-processing and noise removal}

In this section we introduce multiplicative noise to the data to simulate
noise that appears in real measurements 
\begin{equation}
u^{\xi }\left( 0,t\right) =u\left( 0,t\right) \left( 1+\text{rand}\left( %
\left[ -\xi ,\xi \right] \right) \right) ,\text{ }u_{x}^{\xi }\left(
0,t\right) =u_{x}\left( 0,t\right) \left( 1+\text{rand}\left( \left[ -\xi
,\xi \right] \right) \right) ,  \label{6.6}
\end{equation}%
where $\text{rand}\left( \left[ -\xi ,\xi \right] \right) $ is a random
variable uniformly distributed in the interval $\left[ -\xi ,\xi \right] $.
In all our tests we set $\xi =0.1$, which corresponds to the $10\%$ noise.
Functions $u(0,t),u_{x}(0,t)$ and their noisy analogs $u^{\xi
}(0,t),u_{x}^{\xi }(0,t)$ are depicted on Figures 2(B),(C).

The developed numerical technique requires the function $w(x,t)\in
B(d,p_{0},p_{1})$, see (\ref{3.3}) and by (\ref{2.161}) functions $%
p_{0}\left( t\right) ,p_{1}\left( t\right) $ are obtained via the
differentiation of the data $f_{0}\left( t\right) $ and $f_{1}\left(
t\right) $. Thus, the noisy data (\ref{6.6}) should be smoothed out by an
appropriate procedure. To do the latter, we use the cubic smoothing spline
interpolation satisfying the following end conditions: 
\begin{equation*}
\begin{aligned} u(0,0) = 0.5, \hspace{0.3em} u_{tt}(0,T) = 0, \hspace{0.3em}
u_x(0,0) = 0, \hspace{0.3em} u_{xtt}(0,T) = 0. \end{aligned}
\end{equation*}%
Next, we differentiate so smoothed functions. Our numerical experience tells
us that this procedure works quite well.\ Similar observations took place in
all above cited works on the convexification.

\subsection{Numerical results}

We have calculated the relative error of the reconstruction on the final
iteration $n=n^{\ast }$ of the minimization procedure: 
\begin{equation*}
\mbox{ \textit{error} }=\Vert a_{n^{\ast }}-a^{\ast }\Vert
_{L^{2}(0,1)}/\Vert a^{\ast }\Vert _{L^{2}(0,1)}
\end{equation*}
where $a_{comp}(x) = a_{n^{\ast }}(x)$ is the computed solution and $a^{\ast
}(x)$ is the true test function.

We have conducted our computations for the following four tests:

\vspace{1em} \textbf{Test 1}. $a(x)=x^{2}\hspace{0.1em}e^{-(2x-1)^{2}}.$ 
\vspace{1em}

\noindent The function of Test 1 is our reference function for which we have
chosen the above listed parameters. We have used the same parameters in the
remaining Tests 2-4.

\vspace{1em} \textbf{Test 2}. $a(x)=10\hspace{0.1em}e^{-100(x-0.5)^{2}}.$

\textbf{Test 3}. $a(x)=2e^{-400(x-0.3)^{2}}+\hspace{0.1em}%
2e^{-200(x-0.5)^{2}}+\hspace{0.1em}2e^{-400(x-0.7)^{2}}.$

\textbf{Test 4.} $a(x)=1-\sin \left(\frac{ \pi (x-0.876)}{ 1+\pi (x-0.876)}%
\right).$ \vspace{1em}

Note that functions on the Figures 3(C),(D) do not attain zero values at $x=1$ as required by condition (\ref{1.2}). Also note that the function $a\left(x\right) $ in Test 4 is not differentiable at $x_{0}=0.876-\pi ^{-1}\approx 0.558,$ and has infinitely many oscillations in the neighborhood of the point $x_{0}$. Nevertheless numerical reconstructions on Figures 3(A),(D) are rather good ones, also, see Table 6.2. Graphs of exact and computed functions $a(x)$ of Tests 1-4 are presented on Figures 3 (A)-(D). Table 6.2
summarizes the results of our computations.

We have used the 12-core Intel(R) Xeon(R) CPU E5-2620 2.40GHz computer. The
average computational time for tests 1-4 was 159.4 seconds with the
parallelization of our code. And it was 1114.3 seconds without the
parallelization. Thus, the parallelization has resulted in about 7 times
faster computations. \vspace{1em}

\noindent \textbf{Table 6.2. Summary of numerical results. Here }$\left\Vert
\cdot\right\Vert _{\infty }$ \textbf{denotes the }$L_{\infty }$\textbf{\
norm.}

\begin{tabular}{|p{0.75cm}|p{0.5cm}|p{1.25cm}|p{2.2cm}|p{2.2cm}|}
\hline
\emph{Test} & $n^{\ast}$ & \emph{Error} & $\Vert
J^h_{\lambda,\beta,\mu}(w_0) \Vert_{\infty}$ & $\Vert
J^h_{\lambda,\beta,\mu}(w_{n^{\ast}}) \Vert_{\infty}$ \\ \hline\hline
\emph{1} & 30 & 0.1628 & 2570 & 2.7465 \\ \hline
\emph{2} & 33 & 0.2907 & 34.42 & 0.22 \\ \hline
\emph{3} & 51 & 0.0804 & 3.12 & 0.0007 \\ \hline
\emph{4} & 41 & 0.3222 & 0.82 & 0.0003 \\ \hline
\end{tabular}
\vspace{1em}

One can see from Table 6.2 that that
the $L_{\infty }-$norm of the functional $J_{\lambda
,\beta ,\mu }^{h}$ decreases by at least the factor of $150$ in all tests. The same was observed for the $L_{\infty }-$norm of the gradient of this functional (not shown in the table).

\begin{figure}[tbp]
\begin{center}
\subfloat[\emph{$a(x) = 10 \hspace{0.1em}
e^{-100(x-0.5)^2}$}]{\includegraphics[width=.48\textwidth]{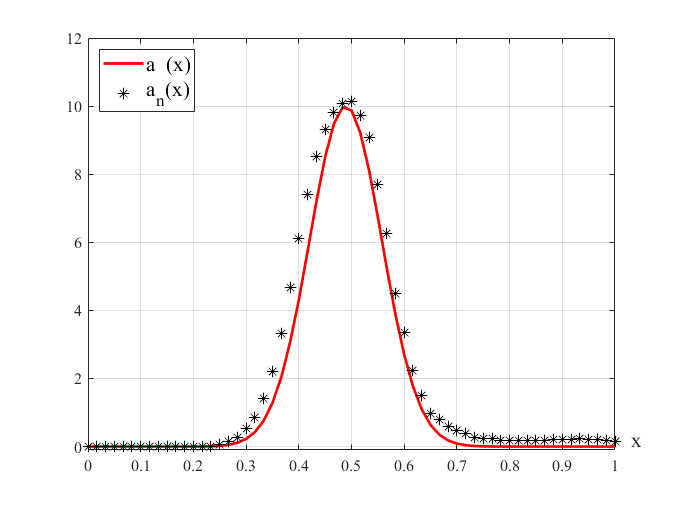}} \quad 
\subfloat[{\emph{$a(x) = 2 ( \hspace{0.1em}e^{-400 (x-0.3)^2} +
\hspace{0.1em} e^{-200 (x-0.5)^2} + \hspace{0.1em}e^{-400
(x-0.7)^2})$}}]{\includegraphics[width =.48\textwidth]{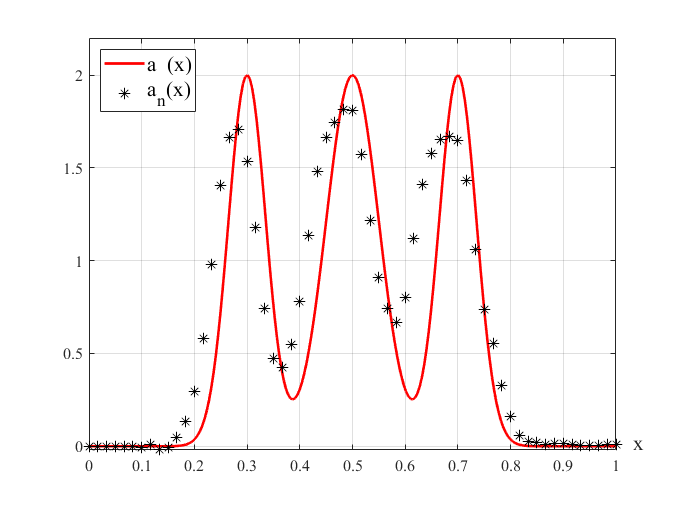}} \\[0pt]
\subfloat[\emph{$a(x) = x^2 \hspace{0.1em}
e^{-(2x-1)^2}$}]{\includegraphics[width =.48\textwidth]{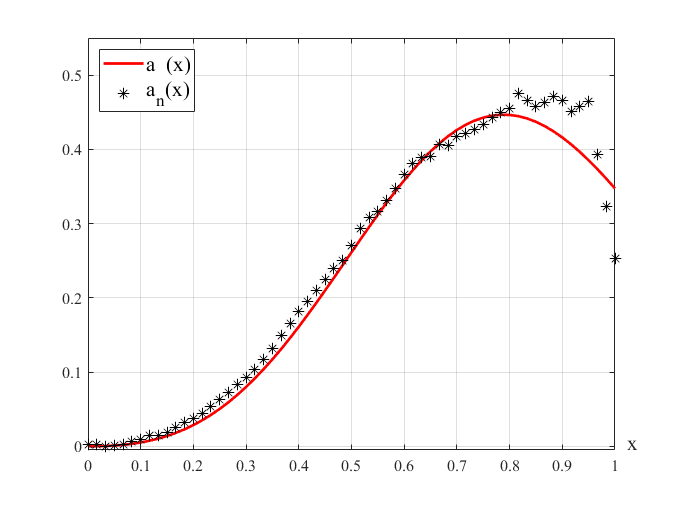}} \quad 
\subfloat[{\emph{$a(x)=1-\sin \left(\frac{ \pi (x-0.876)}{ 1+\pi
(x-0.876)}\right)$}}]{\includegraphics[width=.48\textwidth]{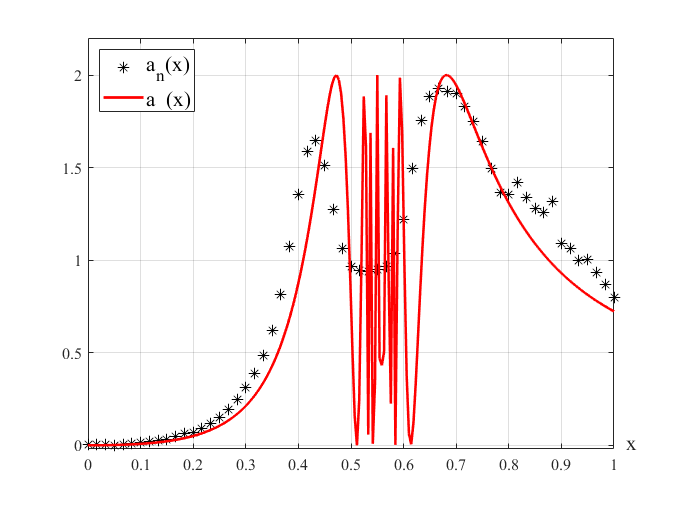}}
\end{center}
\caption{\emph{Numerical reconstructions (the black marked dots) of
functions $a(x)$ (the solid lines). Noise level $\protect\xi =0.1$.}}
\label{fig3}
\end{figure}

We now test some values of the parameters $\lambda $ and $\alpha $ which are
different from ones in (\ref{9}).
Below we work only with the noiseless data and with the function $a\left(
x\right) $ which was used in Test 1. The parameter $\beta $ below is the
same as in (\ref{9}).

First, we test values of $\lambda $ which are larger and smaller than $%
\lambda =2$ in (\ref{9}). Figure 4(A) shows our result for $\lambda =5$ and $%
\lambda =1.$ It is clear from Figure 4(A) that a larger value of $\lambda =5$ provides basically the same result as the one on Figure 3(A), and both are close to the correct solution. On the other hand, the result deteriorates
for a smaller value $\lambda =1.$ Next, Figure 4(B) displays our result for
the limiting case of $\lambda =0,$ i.e. when the Carleman Weight Function is
absent in functional (\ref{3.4}). In this case the gradient descent method
diverges, see Figure 2(D). Thus, we stop iterations after $n=10$ steps. A
significant deterioration of the result of Figure 4(B), as compared with
Figures 3(A) and 4(A), is evident. Therefore, the presence of the CWF in (%
\ref{3.4}) is important.

\begin{figure}[tbp]
\begin{center}
\subfloat[\emph{Exact $a\left( x\right) $ (solid line), $\lambda =1$ (dashed line), $%
\lambda =2$ (dotted line) and $\lambda =5$ (star line).}]{\includegraphics[width=.48\textwidth]{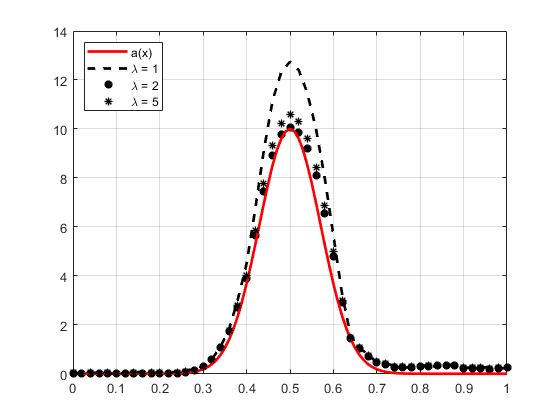}}
\quad \subfloat[{\emph{$\lambda =0$ (dotted line), solid line depicts the
true solution.}}]{\includegraphics[width =.48\textwidth]{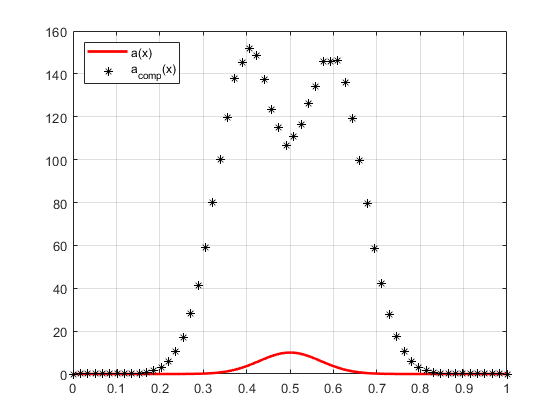}}
\end{center}
\caption{\emph{Limiting testing of different values of the parameter $%
\protect\lambda $ for the test function of Test 1, see comments in the text.
The data are noiseless.}}
\label{fig4}
\end{figure}

\begin{figure}[tbp]
\begin{center}
\subfloat[\emph{$a_{comp}(x)$ for Test 1 with $\alpha = 0.2$ (dashed line)
and $\alpha = 0.5$ (dotted 
line).}]{\includegraphics[width=.6\textwidth]{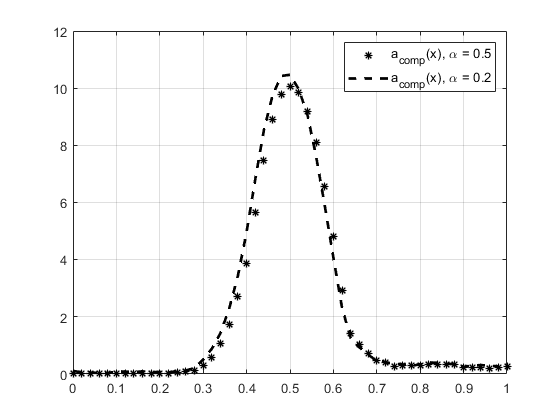}}
\end{center}
\caption{\emph{Testing of different values of the parameter $\protect\alpha ,
$ see comments in the text. Solid line is the correct function of Test 1.
The data are noiseless.}}
\label{fig5}
\end{figure}

The parameter $\alpha $ is chosen in the interval $(0,0.5)$. Figure 5 shows
our results for two values of $\alpha =0.2$ and $\alpha =0.5.$ Here, $%
\lambda =2$, as in (\ref{9}). One can see that both results are almost the
same. A similar behavior was observed for $\alpha =0.3$ and $\alpha =0.4.$
This shows a good stability of our method with respect to the value of $%
\alpha .$ We note that we have chosen the limiting value $\alpha =0.5$ in
our above tests in order to demonstrate that our method is robust with
respect to the choice of $\alpha $ even if the limiting value of this parameter is chosen.

\bibliography{mybib}
\bibliographystyle{plain} 

\medskip 
Received xxxx 20xx; revised xxxx 20xx. \medskip

\end{document}